\documentclass[11pt,a4paper]{article}
\usepackage{amssymb, amsmath, latexsym, proof}

\usepackage{xcolor,soul}

\newtheorem{theorem}{Theorem}
\newtheorem{lemma}{Lemma}[section]

\newtheorem{corollary}[lemma]{Corollary}
\newtheorem{proposition}[lemma]{Proposition}

\newtheorem{definition}{Definition}
\newtheorem{example}[lemma]{Example}
\newtheorem{remark}[lemma]{Remark}

\newcommand{\bl}{\begin{lemma}}
\newcommand{\el}{\end{lemma}}
\newcommand{\bt}{\begin{theorem}}
\newcommand{\et}{\end{theorem}}
\newcommand{\bcor}{\begin{corollary}}
\newcommand{\ecor}{\end{corollary}}
\newcommand{\bp}{\proof{.}}
\newcommand{\ep}{\eop}
\newcommand{\bpr}{\begin{proposition}}
\newcommand{\epr}{\end{proposition}}
\newcommand{\brem}{\begin{remark} \em}
\newcommand{\erem}{\end{remark}}
\newcommand{\bd}{\begin{definition} \em}
\newcommand{\ed}{\end{definition}}
\newcommand{\bex}{\begin{example} \em
}
\newcommand{\eex}{\end{example}}
\newcommand{\beq}{\begin{equation} }
\newcommand{\eeq}{\end{equation}}

\newcommand{\bi}{\begin{itemize}
  }
\newcommand{\ei}{\end{itemize}}
\newcommand{\ben}{\begin{enumerate} }
\newcommand{\een}{\end{enumerate} }

\newcommand{\refeq}[1]{(\ref{#1})}
\newenvironment{enumr}{

\begin{enumerate}     }{\end{enumerate}

}

\newcommand{\benr}{\begin{enumr}
  }
\newcommand{\eenr}{
\end{enumr}}

\newcommand{\ignore}[1]{}

\newcommand{\al}[1]{\forall #1\:}
\newcommand{\ex}[1]{\exists #1\:}

\newlength{\hilflh}

\newcommand{\naturals}{\mathbb{N}}

\renewcommand{\emptyset}{\varnothing}

\newcommand{\cL}{{\mathcal L}}

\newcommand{\cF}{{\mathcal F}}

\newcommand{\fG}{\mathfrak{G}}

\newcommand{\ga}{\alpha}
\newcommand{\gb}{\beta}

\renewcommand{\ge}{\varepsilon}
\newcommand{\gl}{\lambda}
\newcommand{\gL}{\Lambda}

\newcommand{\gy}{\gamma}
\newcommand{\gw}{\omega}
\newcommand{\gS}{\Sigma}

\renewcommand{\phi}{\varphi}

\newcommand{\imp}{\rightarrow}
\newcommand{\eqv}{\leftrightarrow}

\newcommand{\Tr}{\mathrm {Tr}}
\newcommand{\LL}{\mbox{\textit{\textbf{PL}}}}

\newcommand{\Glp}{\mathrm{GLP}}

\newcommand{\Con}{\mathrm{Con}}

\newcommand{\Prf}{\mathrm{Prf}}

\newcommand{\PRA}{\mathsf{PRA}}

\newcommand{\RFN}{\mathsf{RFN}}

\newcommand{\PA}{\mathsf{PA}}
\newcommand{\EA}{\mathrm{EA}}

\newcommand{\ACA}{\mathsf{ACA}}

\newcommand{\gn}[1]{\ulcorner #1 \urcorner}

\newcommand{\ord}[2]{|#2|_{\Pi_1^0}}

\newcommand{\la}{\langle}
\newcommand{\ra}{\rangle}

\renewcommand{\models}{\vDash}      

\newcommand{\On}{\mathrm{On}}
\newcommand{\Cr}{\mathrm{Cr}}
\newcommand{\en}{\mathrm{en}}

\newcommand{\Lim}{\mathrm{Lim}}

\newcommand{\iffdef}{\stackrel{\text{def}}{\iff}}

\newcommand{\nat}{\naturals}

\renewcommand{\leq}{\leqslant}
\renewcommand{\geq}{\geqslant}

\newcommand{\Rc}{\mathrm{RC}}

\newcommand{\dom}{\text{dom}}

\newcommand{\eop}{$\Box$ \protect\par \addvspace{\topsep}}
\newcommand{\proof}[1]{\protect\par\addvspace{\topsep}\noindent {\bf Proof#1}}

\newcommand{\UTB}{\mathsf{UTB}}
\newcommand{\CT}{\mathsf{CT}}
\newcommand{\TT}{\mathsf{T}}
\newcommand{\num}{\underline}
\newcommand{\Set}{\mathrm{Set}}

\newcommand{\tRFN}[2]{{#1}\textrm{-}\RFN({#2})}
\newcommand{\tCA}{\textrm{-CA}}
\newcommand{\lh}{\mathrm{lh}}
\newcommand{\bit}{\mathrm{bit}}

\renewcommand{\LL}{\mathcal{L}}
\newcommand{\tR}{\mathsf{R}}

\newcommand{\bgvee}{\textstyle\bigvee}

\newcommand{\Wo}{\mathbb{W}}
\newcommand{\Fo}{\mathbb{F}}

\newcommand{\Rcl}{{\Rc_\gL}}

\newcommand{\IB}{\mathsf{IB}}
\newcommand{\ATR}{\mathsf{ATR}}

\renewcommand{\ord}{\mathrm{ord}}

\begin{document}

\title{Reflection algebras and conservation results for theories of iterated truth\thanks{Research financed by a grant of the Russian Science Foundation (project No.~16--11--10252). } }
\author{Lev D. Beklemishev\footnote{The corresponding author,  e-mail: \texttt{bekl@mi-ras.ru}.} \  and Fedor N. Pakhomov\\ 
Steklov Mathematical Institute of Russian Academy of Sciences\\ Gubkina str. 8, 119991 Moscow, Russia}
\maketitle 
\begin{abstract}
We consider extensions of the language of Peano arithmetic by transfinitely  iterated truth definitions satisfying uniform Tarskian biconditionals. Without further axioms, such theories are known to be conservative extensions of the original system of arithmetic. Much stronger systems, however, are obtained by adding either induction axioms or reflection axioms on top of them. Theories of this kind can interpret some well-known predicatively reducible fragments of second-order arithmetic such as iterated arithmetical comprehension.

We obtain sharp results on the proof-theoretic strength of these systems using methods of provability logic. Reflection principles naturally define unary operators acting on the semilattice of axiomatizable extensions of our basic theory of iterated truth. The substructure generated by the top element of this algebra provides a canonical ordinal notation system for the class of theories under investigation. 
\ignore{The feature that canonical ordinal notations are not given externally but are extracted from the given structure of reflection principles 
is characteristic for the present approach to ordinal analysis.} 

Using these notations we obtain conservativity relationships for iterated reflection principles of different logical complexity levels corresponding to the levels of the hyperarithmetical hierarchy, i.e., the analogs of Schmerl's formulas. 
These relationships, in turn, provide proof-theoretic analysis of our systems and of some related predicatively reducible theories. In particular, we uniformly calculate the ordinals characterizing the standard measures of their proof-theoretic strength, such as provable well-orderings, classes of provably recursive functions, and $\Pi_1^0$-ordinals.
\end{abstract}

\section{Introduction}

This paper continues the line of research that started with \cite{Bek04}. The main idea of this approach is to study questions of proof-theoretic analysis of theories from the point of view of provability logic methods. Using such methods leads one to consider algebraic structures that occupy, in a sense, an intermediate position between the highly complex world of formal theories (such as systems of first- or second-order arithmetic and set theory) and the simpler and constructive world of natural ordinal notation systems. 

This approach has previously been applied to Peano arithmetic, its fragments and modest extensions. In this paper we make the next necessary step and consider from this point of view theories of predicative strength. Theories of this kind emerged in the works of Solomon Feferman and Kurt Sch\"utte in the 1960s who explicated the informal notion of predicative proof in terms of certain systems of ramified analysis and isolated the proof-theoretic ordinal $\Gamma_0$ as a bound to transfinite induction provable in such  systems \cite{Fef64,Schu}. Since that time, systems reducible to predicative ones in the sense of Feferman and Sch\"utte have been studied quite extensively.

Even though the emphasis in the work in proof theory later shifted to much   stronger impredicative theories, theories of predicative strength and bounds thereof remain an important landmark. Thus, a system $\ATR_0$ of proof-theoretic strength exactly $\Gamma_0$ was isolated by Harvey Friedman in pursuing the reverse mathematics program. Modulo a weak second order theory, that system turned out to be equivalent to some well-known theorems of ordinary mathematics.

Feferman with his collaborators returned to the analysis of predicativity several times during his long career, formulating ever more convincing and simpler to state predicative systems (see e.g.~\cite{Fef91,FefStr00}). More recently, Nick Weaver~\cite{Wea07} disputed the fact that systems presented by Feferman and Sch\"utte adequately represent the informal notion of predicativity. We are sympathetic with the doubts of this kind. However, we also believe that the class of systems isolated by Feferman and Sch\"utte is important in its own right irrespectively of the association of one or another informal notion of predicative proof with it. Henceforth, we will be using the term predicative in the sense of Feferman--Sch\"utte. 

The goal of this paper is to provide an alternative approach to proof-theoretic analysis of the class of theories of proof-theoretic strength below $\Gamma_0$. The main feature of this analysis is the use of hierarchies of reflection principles and conservation results between them instead of the use of infinitary calculi. Under this approach, the ordinal notation systems emerge in the guise of weak propositional (strictly positive) logics, i.e., fragments of provability logic representing the algebras of reflection principles. The close relation between the objects representing the ordinals and the objects representing the theories makes for us the reductions from one to the other rather simple. 

Provalibility logics with transfinitely many modalities, such as $\Glp_\Lambda$ and their strictly positive counterparts $\Rcl$, have been studied already for some time, most notably in the series of works by Joost Joosten and David Fern\'andez-Duque, see~\cite{Bek05a,JF13,JF12-aiml,JF14, BFJ14,Joo15,JooRey16,JF18}.  However, our study seems to be the first one where this machinery is developed enough to yield proof-theoretic analysis results for theories significantly stronger than Peano arithmetic. The main novelty of this paper from a technical point of view are basic conservation theorems (Theorems \ref{th:llpicons} and \ref{redt2}) generalizing the so-called reduction property in arithmetic. The machinery of reflection algebras is then used to generalize these results to transfinitely iterated reflection principles and higher levels of logical complexity. 

Another advantage of our approach is that it allows one to obtain the main results associated with proof-theoretic analysis of theories --- such as consistency proofs, determination of the classes of provably total computable functions, characterization of provable well-orderings --- all at once. This is based on a general formula relating hierarchies of reflection principles of different strength. Formulas of this kind appeared in the work of Ulf Schmerl, first for arithmetic \cite{Schm} and then for the ramified analysis~\cite{Sch82}. In this paper we are, in a sense, redoing the work of Schmerl, but in a different way and, to simplify matters, for a different class of systems (which are first order rather than ramified second order). 

The reflection principles for theories of truth already have been investigated (see \cite[Chapter~XII]{Cie17}). Graham Leigh~\cite{Leigh16} have studied  reflection principles over theories of iterated truth predicates. He characterizes arithmetical consequences of such theories in terms of transfinite induction using the tools of infinitary proof theory. Our basic theory of uniform Tarski biconditionals is rather similar to the one of~\cite{Leigh16}, however we consider finer hierarchies of reflection principles and use different methods.  

Our treatment is essentially self-contained and includes sometimes careful proofs of results that are more or less well-known. The only proof-theoretic tool that we are using is the standard cut-elimination for predicate logic in the form of Tait calculus. We also rely on the method of arithmetization along with the construction of truth definitions. 

The basic systems for which we state our results are formulated in the first-order language of Peano arithmetic expanded by a series of new unary predicate letters representing iterated truth definitions. These predicates are postulated to satisfy rather weak axioms for truth --- the so-called uniform Tarski biconditionals.\footnote{From the point of view of the theory of truth it is more interesting and common to consider stronger axioms for truth predicates, e.g., the compositional axioms $\CT$. Theories based on compositional axioms for truth can also be considered in our context, and in fact they were considered by Lev Beklemishev and Evgeny Dashkov at an earlier stage of this project. However, this choice leads to undesirable complications, both on the level of proof theory and on the level of modal logic, which we prefer to avoid.} The use of weak axioms for truth and delegating the main power of the system to reflection principles added on top of them is one of the key technical ideas of this paper. It is interesting to remark that in this way we can recover some of the results about the compositional truth axioms, e.g., a theorem due to Henryk Kotlarski on the strength of $\Delta_0(\TT)$-induction schema over $\CT$~\cite{Kot86}.    
The plan of this paper is as follows. In Sections \ref{red1} and \ref{red2} we prove two key conservation results for the theory of uniform Tarski biconditionals extended by restricted uniform reflection principles. The first of these results gives a conservation of $\Pi_1(\TT)$-reflection principle over the arithmetical uniform reflection principle, for all arithmetical sentences. The second result, for each $n>0$, is a conservation of $\Pi_{n+1}(\TT)$-reflection principle over $\gw$-times iterated $\Pi_n(\TT)$-reflection principle for $\Pi_n(\TT)$-sentences. In Section \ref{it-tr} we introduce theories of iterated truth predicates and the corresponding reflection principles. We reinterpret the results of Sections \ref{red1} and \ref{red2} for these theories as two different reduction formulas. Section \ref{rcl} considers a strictly positive logic with transfinitely many modalities, $\Rc_\Lambda$, which is a straightforward generalization of the reflection calculus $\Rc$. We recapitulate the results of \cite{Bek05a} calculating the order types of the ordering of $\Rc_\Lambda$ formulas. In Section~\ref{ord-an} we capitalize on the previous work and obtain Schmerl-type conservation results for iterated truth theories by using $\Rc_\Lambda$. In Section \ref{anref} we outline the ideas of ordinal analysis of iterated reflection and the use of reflection algebras for this purpose. In Section \ref{sec-ord} we apply these results to obtain conservation theorems and ordinal analysis of some standard systems of second-order arithmetic of predicative strength.

\paragraph{Acknowledgements.} This paper has a long prehistory. It was planned by the first author as a second part of \cite{Bek04} and a continuation of \cite{Bek05a}, however the work turned out to be more complicated than expected and experienced various delays. Later other people, most notably Joost Joosten and David Fern\'andez, joined the efforts in the study of transfinite extensions of provability logic and their proof-theoretic interpretations. Evgeny Dashkov contributed by his work \cite{Das12en} on strictly positive provability logics that was motivated by the needs to adapt the provability algebras framework to stronger reflection principles. The first author and Dashkov subsequently worked on theories of iterated truth predicates based on Tarski compositional axioms, but the paper was never completed. The present paper also owes a lot to Evgeny Kolmakov whose active interest revitalized this work and who contributed to many discussions of the current approach. To all the people mentioned above we express our sincere gratitude. 

F. Pakhomov was
selected as one of the Young Russian Mathematics award winners, and he would like to thank its sponsors
and jury; his work on this particular paper was funded from another source.

\section{Preliminaries: languages and truth theories}
As our basic system of arithmetic we take \emph{Elementary Arithmetic} $\EA$, also known as EFA or $I\Delta_0(\exp)$, in any of its standard formulations (see \cite{HP,Bek05}). The language of $\EA$ has symbols for successor, addition, multiplication, exponentiation functions and the order relation. The axioms of $\EA$ comprise basic defining equations for all these symbols, as well as the induction axioms for bounded ($\Delta_0$) formulas. We allow $\exp$ in quatifier bounding terms in the definition of the class $\Delta_0$. $\EA^+$ denotes the extension of $\EA$ by an axiom asserting the totality of superexponentiation function, also known as $I\Delta_0+\text{Supexp}$ (see~\cite{HP}). 

Let $\LL$ be a first order language extending that of $\EA$ with at most countably many fresh predicate symbols. We assume fixed an elementary G\"odel numbering of $\LL$.

Let $\Delta^{\LL}_0$ denote the class of formulas obtained from atomic $\LL$-formulas and formulas of the form $\TT(t)$, where $t$ is a term, by Boolean connectives and bounded quantifiers. The classes $\Pi_n^\LL$ and $\Sigma_n^\LL$ are defined from $\Delta^{\LL}_0$ in the usual way:  $\Pi_0^\LL=\Sigma_0^\LL=\Delta_0^\LL$, $\Pi_{n+1}^\LL=\{\al{\vec x}\phi(\vec x):\phi\in\gS^\LL_n\}$, and $\gS_{n+1}^\LL=\{\ex{\vec x}\phi(\vec x):\phi\in\Pi^\LL_n\}$.
$I\Gamma$ denotes the induction schema restricted to $\Gamma$-formulas. 
$\EA^\LL$ denotes the extension of $\EA$ by $\Delta_0^\LL$-induction schema.

The extension of $\LL$ by a unary predicate symbol $\TT$ is denoted $\LL(\TT)$. For typographical reasons the classes $\Delta_0^{\LL(\TT)}$, $\Pi_n^{\LL(\TT)}$, $\Sigma_n^{\LL(\TT)}$ will also be denoted 
$\Delta_0^{\LL}(\TT)$, $\Pi_n^{\LL}(\TT)$, $\Sigma_n^{\LL}(\TT)$, respectively.

We consider an $\LL(\TT)$-theory $\UTB_\LL$ axiomatized by the following schemata:

\bi 
\item[U1:] $\forall \vec x \left( \phi(\vec x) \eqv \TT(\gn{\phi(\num{\vec{x}})}) \right),$ for all $\LL$-formulas $\phi(\vec x)$;
\item[U2:] $\neg \TT(\num{n})$, for all $n$ such that $n$ is not a G\"odel number of an $\LL$-sentence.
\ei 
Here, $\gn{\phi(\num{\vec x})}$ is an elementarily definable term representing the function mapping $n_1,\dots,n_k$ to the G\"odel number of $\phi(\num{n}_1,\dots,\num{n}_k)$, provided $\vec x=x_1,\dots,x_k$.

The theory $\UTB_\LL$ plays a central role in this paper. We remark that it can be given a natural $\Pi_2^\LL$-axiomatization by taking U2 together with the following schemata, for all $\LL$-formulas $\phi$, $\psi$:

\bi 
\item[C1:] $\forall \vec x \left( \phi(\vec x) \eqv \TT(\gn{\phi(\num{\vec{x}})}) \right),
$ if $\phi(\vec x)$ is atomic;
\item[C2:] $\al{\vec x}(\TT(\gn{\phi(\vec{\num x})\land\psi(\vec{\num x})})\eqv \TT(\gn{\phi(\vec{\num x})})\land\TT(\gn{\psi(\vec{\num x})}))$;
\item[C3:] $\al{\vec x}(\TT(\gn{\neg\phi(\vec{\num x})})\eqv \neg\TT(\gn{\phi(\vec{\num x})}))$;
\item[C4:] $\al{\vec x}(\TT(\gn{\al{y}\phi(y,\vec{\num x})})\eqv \al{y}\TT(\gn{\phi(\vec{\num x})}))$.
\ei 

If Axioms C2--C4 are stated as `global' axioms with a universal quantifier over G\"odel numbers of $\LL$-formulas $\phi,\psi$, the corresponding theory is known as \emph{Tarski compositional axioms for truth} and is denoted $\CT_\cL$ in this paper.\footnote{It is denoted $\CT_\cL^-$ in some treatments to stress the absence of induction axioms. We do not presuppose any induction axioms in $\CT_\cL$.} 

\ignore{
Let $\Delta^{\LL}_0(\TT)$ denote the class of formulas obtained from \emph{all} $\LL$-formulas and atomic formulas of the form $\TT(t)$ by Boolean connectives and bounded quantifiers. The classes $\Pi_n^\LL(\TT)$ and $\Sigma_n^\LL(\TT)$ are defined from $\Delta^{\LL}_0(\TT)$ in the usual way:  $\Pi_0^\LL(\TT)=\Sigma_0^\LL(\TT)=\Delta_0^\LL(\TT)$, $\Pi_{n+1}^\LL(\TT)=\{\al{\vec x}\phi(\vec x):\phi\in\gS^\LL_n(\TT)\}$, and $\gS_{n+1}^\LL(\TT)=\{\ex{\vec x}\phi(\vec x):\phi\in\Pi^\LL_n(\TT)\}$.
}

The following lemma is well-known and straightforward, it is verifiable in $\EA$ by building a local interpretation. Let $S$ be any $\LL$-theory containing $\EA$.

\bl $S + \UTB_\LL$ is conservative over $S$ for $\LL$-formulas. \el

By a \emph{G\"odelian theory} we mean a theory, in a language as above, whose set of axioms comes equipped with an elementary ($\Delta_0(\exp)$) formula defining its set of G\"odel numbers in the standard model of arithmetic. With every such theory $S$ we associate its $\Sigma_1$ provability predicate $\Box_S(x)$ in a standard way~\cite{Fef60}. We say that \emph{$S_1$ $U$-provably contains $S_2$} if 
$$U\vdash \al{x}(\Box_{S_2}(x)\to \Box_{S_1}(x)).$$ 

Suppose $S$ is G\"odelian and $\Gamma$ is a set of formulas in the language of $S$. By $\tRFN{\Gamma}{S}$ we denote the \emph{uniform reflection principle for $\Gamma$-formulas}, that is, the schema 
$$\al{x}(\Box_S\phi(\num{x})\to \phi(x)), \qquad \phi\in \Gamma.$$ In particular, $\tRFN{\LL}{S}$ is the uniform reflection principle for all $\LL$-formulas. $\tRFN{\TT}{S}$ is the following $\TT$-reflection principle:
$$
\forall \phi\in\LL\, \left(\Box_S \phi \imp \TT(\phi)\right).
$$

We note two obvious lemmas.
\begin{lemma} \label{utb1}
Suppose $S$ $\EA$-provably contains $\UTB_\LL$.
\benr
\item $\EA \vdash \forall \psi\in\LL\, \Box_{\UTB_\LL} \left(\psi \eqv \TT(\psi)\right)$;
\item $\tRFN{\TT}{S}$ is equivalent to $\forall \phi\in\LL\, \left(\Box_S \TT(\phi) \imp \TT(\phi)\right)$ over $\EA$;
\item     $\EA+\UTB_\LL+\tRFN{\TT}{S}$ contains $\tRFN{\LL}{S}$;
\item $\EA+\tRFN{\Delta^\LL_0(\TT)}{S}$ contains $\tRFN{\TT}{S}$;
\item $\EA+\tRFN{\Pi_2^\LL(\TT)}{S}$ proves the compositional axioms $\CT_\cL$.
\eenr 
\end{lemma}

\bp\ Concerning the proof of (v), we remark that provably in $\EA$, for all $\cL$-formulas $\phi,\psi$, conditions C2--C4 are provable in $S$ (since $S$ contains $\UTB$). Since C2--C4 are at most $\Pi_2^\cL$, we infer $\CT_\cL$ by applying reflection for $S$. \ep 

Let $\LL$ be a language with or without $\TT$. 
\bl \label{pi-ref}
\benr
\item $\tRFN{\Pi_1^\LL}{S}$ is equivalent to $\tRFN{\Delta_0^\LL}{S}$ over $\EA$;
\item $\tRFN{\Pi_{n+1}^\LL}{S}$ is equivalent to $\tRFN{\Sigma_n^\LL}{S}$ over $\EA$.
\eenr 
\el

Next we recall that reflection implies induction, as noted by Kreisel and L\'evy~\cite{KrL}.
\bl  \label{ref-ind} If $S$ is an $\LL$-theory provably containing $\EA$, then
\benr
\item
$\EA + \tRFN{\Delta_0^\LL}{S}\vdash I\Delta_0^\LL$;
\item $\EA+\tRFN{\Pi_{n+2}^\LL}{S}\vdash I\Pi_{n}^\LL$, for each $n\geq 1$;
\item $\EA + \tRFN{\LL}{S}\vdash I\LL$.
\eenr
\el

\bp\ (i) The usual argument goes, for any $\Delta_0^\LL$-formula $\phi$, as follows. First, show that
$$\EA\vdash \al{n}\Box_S (\phi(0)\land \al{x\leq \num{n}}(\phi(x)\to \phi(x+1))\to \phi(\num{n})).$$ Then, using $\tRFN{\Delta_0^\LL}{S}$ we infer
$$\al{n}(\phi(0)\land \al{x\leq n}(\phi(x)\to \phi(x+1))\to \phi(n)),$$ which implies the standard instance of the induction axiom for $\phi$.

The proof of Claims (ii) and (iii) is similar.
\ep

\brem 
For the language of arithmetic, by a well-known result of D.~Leivant~\cite{Lei83}, the theory $I\Pi_n$ is equivalent to $\EA+\tRFN{\Pi_{n+2}}{\EA}$ for each $n\geq 1$. This implies the earlier result by Kreisel and L\'evy that $\EA+\RFN(\EA)\equiv \PA$. 

However, for the language with a truth predicate, full reflection over $\EA+\UTB_\cL$ is stronger than full induction. One can prove by a simple model-theoretic argument that the theory $\EA+\UTB_\cL+I\cL(\TT)$ is a conservative extension of $\PA$. On the other hand, the theory $\EA+\tRFN{\cL(\TT)}{\EA+\UTB_\cL}$ is not: By Lemma~\ref{utb1} (v) it proves the compositional axioms $\CT_\cL$, which already with the weaker $\TT$-reflection imply arithmetical reflection over $\PA$  (see Kotlarski \cite{Kot86}). 
\erem

\ignore{Whenever a language $\LL$ with or without $\TT$ is fixed, we will often abbreviate $\tRFN{\Pi_m^\LL}{S}$ by $\tR_m(S)$. }

\section{Conservativity of $\Pi_1(\TT)$-reflection} \label{red1}

\begin{theorem}\label{th:llpicons}
 Suppose $S$ is a theory in a language extending $\LL(\TT)$ and provably containing $\EA + \UTB_\LL$. Then
 $\EA+\UTB_\LL +  \tRFN{\Pi_1^\LL(\TT)}{S}$ is a conservative extension of $\EA + \tRFN{\LL}{S}$ for $\LL$-formulas.
\end{theorem}

\brem $\EA + \tRFN{\LL}{S}$ contains the $\LL$-fragment of $S$, because if $\phi\in\LL$ and $S\vdash \phi$ then $\EA\vdash \Box_S \phi$ and hence $\EA + \tRFN{\LL}{S}\vdash \phi$.
\erem

\bp\ For each finite fragment $\cF$ of the signature of $\LL$ and a finite set $\Gamma \subseteq \UTB_\cF$ we will construct a non-relativizing and $\cF$-preserving interpretation of the theory
$\EA+\Gamma + \tRFN{\Delta_0^\cF(\TT)}{S}$ in $\EA + \tRFN{\cF}{S}$. By Lemma \ref{pi-ref}(i) this yields a suitable local interpretation of $\EA+\UTB_\LL + \tRFN{\Pi_1^\LL(\TT)}{S}$ which implies  $\LL$-conservativity by compactness.

Axioms of $\Gamma$ are either of the form $\neg \TT(\num n)$ or of the form
$$
\forall x_1 \dots \forall x_n\:(\psi(x_1, \dots, x_n) \eqv \TT(\gn{\psi(\num{x}_1, \dots, \num{x}_n)})),
$$
for some $\psi\in\cF$. Fix an $m < \omega$ such that, for each axiom of $\Gamma$ and the corresponding formula $\psi$, both $\psi$ and $\neg \psi$ are in $\Pi^\cF_m$. It is well-known that there is a truth definition $\Tr_{\Pi_m^\cF}$ in $\EA^\cF$ for $\Pi_m^\cF$-formulas (see Theorem~\ref{tr-pi-m} in the Appendix, where a somewhat sharper result needed later is spelled out).  In the following proof we read $\tR_m(S)$ as $\tRFN{\Pi_{m+1}^\cF}{S}$.   By Lemma~\ref{ref-ind}, $\EA+\tR_0(S)$ contains $\EA^\cF$.

In Lemma \ref{theta} we construct an $\cF$-formula $\theta(x)$ that will serve as the interpretation of $\TT(x)$. Here and below $\dot{\land}, \dot{\neg}, \dot{\bigwedge}$ stand for elementary terms defining syntactical operations on the G\"odel numbers of formulas, as in Feferman \cite{Fef60}. 

\begin{lemma} \label{theta}
There is a $\Pi_{m+2}^\cF$-formula $\theta(x)$ such that the following properties hold provably in $\EA + \tRFN{\cF}{S}$: {\rm
\begin{enumerate}
\item[T0.] $\al{\phi}(\theta(\phi)\to \phi\in\cF)$\quad ($\theta$ defines a set of $\cF$-sentences);
\item[T1.] $\forall \phi\in\cF\, (\Box_S(\phi) \imp \theta(\phi))$ \quad ($\theta$ contains the $\cF$-fragment of $S$);

\item[T2.] $\forall \phi\in\cF\, (\theta(\dot{\neg} \phi) \eqv \neg \theta(\phi))$ and $\forall \phi, \psi\in\cF\, (\theta(\phi \mathop{\dot{\land}} \psi) \eqv \theta(\phi) \land \theta(\psi))$ \\
($\theta$ commutes with propositional connectives);

\item[T3.] $\forall \phi \left(\Tr_{\Pi_m^\cF}(\phi) \imp \theta(\phi)\right)$ \quad ($\theta$ contains all true $\Pi_m^\cF$-sentences).
\end{enumerate}}
\end{lemma}

Before proving this lemma, let us show that these properties are sufficient to prove the translations of all the axioms of $\Gamma + \tRFN{\Delta_0^\cF(\TT)}{S}$ in $\EA + \tRFN{\cF}{S}$.\footnote{With some care, one can check that T0--T3 can, in fact, be verified in $\EA+\tR_{m+3}(S)$.}

If an axiom of $\Gamma$ has the form $\neg \TT(\num n)$, where $n$ is not the G.n.\ of an $\LL$-sentence, then $\neg\theta(\num n)$ is provable by T0. The other axioms of $\Gamma$ translate into $\psi(\vec x) \eqv \theta(\gn{\psi(\num{\vec x})})$. Note that
the implication $\psi(\vec x) \imp \theta(\gn{\psi(\num{\vec x})})$ follows from T3.
Applying T3 to $\neg \psi(\vec x)$ we obtain $\theta(\gn{\neg \psi(\num{\vec x})})$ which yields $\neg\theta(\gn{\psi(\num{\vec x})})$ by T2.

To deal with bounded quantifiers we need a definition and two lemmas.
For each $\Delta_0^\cF(\TT)$-formula $\phi(\vec x)$ we specify an elementarily definable term $\phi^*(\vec x)$ as follows:
\begin{enumerate}
\item $\phi(\vec x)^*=\gn{\phi(\num{\vec x})}$, if $\phi(\vec x)$ is an atomic $\cF$-formula;
\item $\TT(t(\vec x))^*= t'(\vec x)$, where $t'(\vec x)$  provably in $\EA$ satisfies
$$t'(\vec x):=\begin{cases} t(\vec x),\quad \text{if $t(\vec x)$ is a G.n.\ of an $\cF$-sentence}, \\ \gn{0=1},\quad \text{otherwise};
\end{cases}
$$
\item $(\phi\land\psi)^*=(\phi^*\dot{\land} \psi^*)$; $(\neg\phi)^*=\dot{\neg}\phi^*$;
\item $(\forall u\leq t\:\phi(u,\vec x))^*=\dot{\bigwedge}_{i\leq t} \phi^*(i,\vec x)$.
\end{enumerate}

For each $\vec n$, $\phi^*(\vec n)$ denotes the G\"odel number of a sentence which is equivalent to $\phi(\num{\vec n})$ in $\UTB_\cF$. Formalizing this fact in $\EA$ yields

\begin{lemma} \label{star1} For each $\Delta_0^\cF(\TT)$-formula $\phi(\vec x)$, $$\EA\vdash \forall \vec x\:\Box_S(\gn{\phi(\num{\vec x})}\mathop{\dot{\eqv}} \phi^*(\vec x)).$$
\end{lemma}

\bp\
Induction on the build-up of $\phi$. If $\phi\in\cF$ and atomic, the claim is trivial. If $\phi$ has the form $\TT(t(\vec x))$ for some term $t$, we reason in  $\EA$ as follows. Given $\vec x$, if $t(\vec x)$ is the G.n.\ of an $\cF$-sentence $\psi$, then $\Box_S(\TT(\psi)\eqv \psi)$ by Lemma \ref{utb1}, and the claim follows. If $t(\vec x)=n$ is not a G.n.\ of an $\cF$-sentence, then $\phi^*(\vec x)=t'(\vec x)=\gn{0=1}$. By the axioms U2 and the equality axioms, the sentence $\TT(t(\num{\vec{x}}))$ is equivalent to $\TT(\num n)$ and refutable in $S$. Hence,  it is equivalent to $0=1$ and to $\phi^*(\vec x)$.

Boolean connectives preserve the equivalence. For the bounded quantifier, we use the fact that
$$\textstyle\EA\vdash \forall\phi\forall n\ \Box_S(\gn{\forall x\leq \num{n}\:\phi(x)}\mathop{\dot{\eqv}} \dot{\bigwedge}_{x\leq n}\:\gn{\phi(\num{x})}).$$
\ep

Let $\phi^\theta$ denote the translation of a formula $\phi\in\cF(\TT)$ under the substitution $\TT(t)/\theta(t)$.

\begin{lemma} \label{star2} For each $\Delta_0^\cF(\TT)$-formula $\phi$ there holds $$\EA + \tRFN{\cF}{S}\vdash \phi^\theta \eqv \theta(\phi^*).$$
\end{lemma}

\bp\
We argue by induction on the build-up of $\phi$. All cases except for the case of a bounded quantifier are easy. Let $\phi= \forall x\leq t\ \psi(x)$. Then $\phi^\theta$ is equivalent to $\forall x\leq t\ \theta (\psi^*(x))$ by the induction hypothesis.
We claim that
$$\textstyle\EA + \tRFN{\cF}{S}\vdash \forall y\:(\forall x\leq y\ \theta(\psi^*)\eqv \theta(\dot{\bigwedge}_{i\leq y} \psi^*(i)).$$
The implication from left to right is proved by induction on $y$ using T2. By Lemma~\ref{ref-ind} the induction is available in $\EA+\tRFN{\cF}{S}$. \ignore{In fact, since the complexity of $\theta$ is $\Sigma^\cF_{m+1}$, one can infer the relevant instance of induction from $\tR_{m+3}(S)$.} The implication from right to left is easier and can be inferred from T1 and T2 without the use of induction.

From this claim we conclude that $\forall x\leq t\ \theta (\psi^*(x))$ is equivalent to
$\theta(\dot{\bigwedge}_{i\leq t} \psi^*(i))$ which is the same as  $\theta(\phi^*)$. This concludes the proof of lemma.
\ep

\begin{lemma}
$\EA+\tRFN{\cF}{S}$ proves the translation of $\tRFN{\Delta_0^\cF(\TT)}{S}$.
\end{lemma}

\bp\   Let $\phi(x)$ be a $\Delta^\cF_0(\TT)$-formula. Reasoning in $\EA+\tRFN{\cF}{S}$ assume $\Box_S\phi(\num{x})$. Then $\Box_S\phi^*(x)$, by Lemma \ref{star1}. By T1, we infer $\theta(\phi^*(x))$, which implies $\phi(x)^\theta$ by Lemma \ref{star2}.
\ep

Now we prove Lemma \ref{theta}. The formula $\theta(x)$ is constructed via the standard (arithmetized) process of completion of $S$, but with $m$-consistency instead of the usual consistency (cf.~\cite[Theorem 4.11]{Fef60}). We note that this construction is simpler than (and, in fact, a part of) the equally well-known construction of a Henkin-completion of $S$. Namely, let $\phi_0, \phi_1, \dots, \phi_n, \dots$ be an enumeration of all $\cF$-sentences, set $T_0 = \emptyset$ and define
$$
T_{n + 1} := \begin{cases}
T_n + \phi_n, \text{ if $S + \mbox{all true $\Pi^{\mathcal{F}}_m$-sentences} + T_n + \phi_n$ is consistent,}\\
T_n, \text{ otherwise.}
\end{cases}
$$
Next we will give the formula $\theta(x)$ such that $\theta(\varphi_x)$ expresses the fact that $\varphi_x\in T_{\omega}=\bigcup\limits_{n\in\omega}T_n$. 

Let $x\in_A z$ express in $\EA$ that the $x$-th bit in the binary expansion of $z$ is 1. Formally we put $\theta(x)$ to be $\exists y \;(x=\varphi_y\land \theta'(y))$, where $\theta'(x)$ is 
$$\exists s,l(x\le l \land x\in_A s\land (\forall i\le l)(i\in_A s\mathrel{\leftrightarrow} \mathsf{R}_m(S+\{\varphi_j\mid j<i,\;j\in_As\} +\varphi_i))).$$

Now let us check that the theory $\mathsf{EA}+\mathcal{F}\mbox{-}\mathsf{RFN}(S)$ proves the conditions T0--T3 for $\theta(x)$. The condition T0 trivially follows from the definition of $\theta(x)$. Recall that by Lemma \ref{ref-ind} the theory $\mathcal{F}\mbox{-}\mathsf{RFN}(S)$ contains full scheme of $\mathcal{F}$-induction. Reasoning in $\mathsf{EA}+\mathcal{F}\mbox{-}\mathsf{RFN}(S)$ we prove by induction on $n$ that $R_m(S+T_n)$ (which is just a formalized version of the assertion $S+ \mbox{all true $\Pi^{\mathcal{F}}_m$-sentences}+T_n$ is consistent). Using this it is easy to derive that $S\subseteq T_{\omega}$, which yields T1. Similarly, $T_{\omega}$ contains all true $\Pi_m^{\mathcal{F}}$ sentences, which gives us T3. Further, for each $\varphi$ exactly one of the following two possibilities holds 1. $\varphi\in T_{\omega}$ or 2. $\lnot\varphi\in T_{\omega}$; this gives us the part of T2 about negation. Finally, we prove that $T_{\omega}$ is closed under first-order deduction, which in particular implies the part of T2 about conjunction.

This completes the proof of Lemma \ref{theta} and of Theorem \ref{th:llpicons}.
\ep

\ignore{
\begin{corollary}\label{cor:rfncons}
$\UTB + \tRFN{\Pi_1(\TT)}{S}$ is a conservative extension of $\EA + \RFN(S)$ for arithmetical sentences.
\end{corollary}
}

\begin{corollary}\label{cor:rfncons}
If $S$ is an $\LL(\TT)$-theory, then 
$$\EA+\UTB_\LL +  \tRFN{\Pi_1^\LL(\TT)}{S + \UTB_\LL}$$ is a conservative extension of $\EA + \tRFN{\LL}{S}$ for $\LL$-formulas.
\end{corollary}

\bp\
Since  $S + \UTB_\LL$ is conservative over $S$, we have $\tRFN{\LL}{S + \UTB_\LL} \equiv \tRFN{\LL}{S}$.
The result now follows by applying the previous theorem to the theory $S + \UTB_\LL$.
\ep

\section{Conservativity of $\Pi_{n+1}(\TT)$-reflection} \label{red2}

The main result of this section is a relativization of the reduction property for arithmetical reflection principles, that is, of Theorem 2 in \cite{Bek99b}. In the present setup the usual argument in arithmetic goes through without any substantial changes. However, we need to take care of some extra details. The most important detail is the need for the reflection schemata to be finitely axiomatizable, for the easy (inclusion) part of the theorem. This part is based on the existence of partial truth definitions and the finiteness of the language and requires $\Delta_0^\LL(\TT)$-induction, whereas the conservation part does not.

In this section $\LL$ can be a language with or without $\TT$. 
Let $S$ be a G\"odelian theory in a language extending $\LL$ and  provably containing $\EA$. We abbreviate by  $\tR_{S,n}(\phi)$ the schema $\tRFN{\Pi_{n+1}^\LL}{S+\phi}$. We also write $\tR_{S,n}$ for $\tR_{S,n}(\top)$. If $\LL$ is finite, $\tR_{S,n}$ will be finitely axiomatizable. However, in general, $\tR_{S,n}(\phi)$ is an elementarily axiomatized (possibly infinite) schema. We read $R_1\vdash R_2$ as: each instance $\phi_2\in R_2$ is provable in $\EA+R_1$; $R_1\land R_2$ denotes $R_1\cup R_2$, and $R_1\lor R_2$ is the set $\{\phi_1\lor\phi_2:\phi_1\in R_1,\phi_2\in R_2\}$.

\bl \label{Rn} For all $n\geq 0$, the following properties hold provably in $\EA$: For all sentences $\phi,\psi\in\cL$, 
\benr
\item If $S\vdash\phi\to \psi$ then $\tR_{S,n}(\phi)\vdash \tR_{S,n}(\psi)$;
\item $\tR_{S,n}(\phi\lor\psi)\vdash \tR_{S,n}(\phi)\lor \tR_{S,n}(\psi)$;
\item $\tR_{S,n}(\phi)\vdash \phi$ if $\phi\in \Pi^\LL_{n+1}$;
\item $\tR_{S,n}(\phi)\vdash \Diamond_{S}\phi$.
\eenr
\el

\bp\ For each item we give an informal argument that can be readily formalized in $\EA$.

(i) To derive $\tR_{S,n}(\psi)$ assume $\Box_{S+\psi} \theta(\num{m})$ with $\theta\in \Pi_{n+1}^\LL$. Since $S\vdash\phi\to\psi$ we have $\Box_{S+\phi} \theta(\num m)$. Hence,  using $\tR_{S,n}(\phi)$ we obtain $\theta(m)$ by reflection. 

(ii) Assume  $\theta \in \tR_{S,n}(\phi)\lor \tR_{S,n}(\psi)$, then $\theta\circeq\theta_1\lor\theta_2$ with $\theta_1\in\tR_{S,n}(\phi)$ and $\theta_2\in\tR_n(\psi)$. For some formulas $\phi_1,\psi_1\in \Pi_n^\LL$ we have $$\text{$\theta_1\circeq \al{x}(\Box_{S+\phi}\phi_1(\num x)\to \phi_1(x))$ and 
$\theta_2\circeq \al{x}(\Box_{S+\psi}\psi_1(\num x)\to \psi_1(x))$}.$$ We observe that $\theta_1\lor\theta_2$ is logically equivalent to \beq \label{disj} \al{x}\al{y}(\Box_{S+\phi}\phi_1(\num x)\land \Box_{S+\psi}\psi_1(\num y) \to (\phi_1(x)\lor \psi_1(y))).\eeq 
Reasoning in $\EA+\tR_{S,n}(\phi\lor\psi)$ we now prove formula \refeq{disj}. Consider any $x,y$ and assume $\Box_{S+\phi}\phi_1(\num x)\land \Box_{S+\psi}\psi_1(\num y)$. Then obviously $\Box_{S+\phi\lor\psi}(\phi_1(\num x)\lor\psi_1(\num y))$. Hence, $\phi_1(x)\lor\psi_1(y)$ by reflection, which proves \refeq{disj}.

(iii) If $\phi\in\Pi_{n+1}^\LL$, then the schema $\tR_{S,n}(\phi)$ contains $\Box_{S+\phi}\phi\to \phi$. Since $\EA\vdash \Box_{S+\phi}\phi$, $\phi$ follows. 

(iv) The instance of $\tR_{S,n}(\phi)$ for $\bot$ implies  $\Diamond_S\phi$.  
\ep

\bt \label{redt2}
Let a theory $U$ be axiomatized over $\EA$ by a set of $\Pi^\LL_{n+2}$-sentences. Then,
$U + \tR_{S,n+1}$ is a $\Pi^\LL_{n+1}$-conservative extension of the closure of $U$ under the rule $\phi/\tR_{S,n}(\phi)$ for $\phi\in\Pi^\LL_{n+1}$.
\et

\brem
Notice that $\UTB_\LL$ is axiomatized over $\EA$ by $\Pi^\LL_2(\TT)$-sentences, hence the result applies to $\LL(\TT)$-theories $U$ containing $\UTB_\LL$. Also notice that the conclusion $\tR_{S,n}(\phi)$ of the rule is, generally, a schema rather than a sentence, which means that one is allowed to infer from $\phi$ any instance of that schema.
\erem

\bp\
As in the arithmetical case the proof goes by considering a cut-free derivation in Tait calculus of the sequent
\beq
\neg U,\ \Sigma,\  \Pi, \label{seqe}
\eeq
where $\Pi$ is a set of $\Pi^\LL_{n+1}$-formulas, $\neg U$ is a finite set of the negations of the axioms of $U$, and $\Sigma$ is a finite set of the negations of the $\tRFN{\Sigma^\LL_{n+1}}{S}$ schema instances. Every such instance has the form
$$\exists y\exists x\:[\Prf_S(y,\gn{\neg\phi(\num x)})\land \phi(x)],$$
for some $\phi(x)\in{\Pi}^\LL_{n+1}$. Let $P_{\phi}(x,y)$ denote the
formula in square brackets. We can also assume the axioms of $U$ to
have the form $\forall x_1\ldots\forall x_m \neg A(x_1,\ldots,x_m),$
for some $\Pi^\LL_{n+1}$-formulas $A(\vec x)$.

By the subformula property, any formula occurring in a derivation of
a sequent $\Gamma$ of the form \refeq{seqe} either (a) is a
$\Pi^\LL_{n+1}$-formula, or (b) has the form
$\Sigma$, $\exists x P_{\phi}(t,x)$ or
$P_{\phi}(t,s),$ for some terms $s,t$, or (c) has the form
$$\exists x_{i+1}\ldots\exists x_m A(t_1,\ldots,t_i,x_{i+1},\ldots,x_m),$$
for some $i<m$ and terms $t_1,\ldots,t_i$. Let $\Gamma^-$ denote the
result of deleting all formulas of types (b) and (c) from $\Gamma$. Let $U'$ denote the closure of $U$ under the rule $\phi/R_{S,n}(\phi)$ with $\phi\in\Pi^\LL_{n+1}$.

\begin{lemma} \label{genz}
If a sequent $\Gamma$ of the form {\em \refeq{seqe}} is cut-free
provable, then $\bgvee\Gamma^-$ is provable in $U'$.
\end{lemma}
\bp\ Induction on the height $d$ of a derivation of
$\Gamma$. It is sufficient to consider the case that a formula of
type (b) or (c) is introduced by the last application of a rule in
$d$; besides, it is sufficient to only consider the formulas
 $P_{\phi}(t,s)$ and $\exists
x_m A(t_1,\ldots,t_{m-1},x_m)$, because in all the other cases the
premise and the conclusion of the rule coincide, after applying the
operation $(\cdot)^-$. The second case is easy (see \cite{Bek99b}).

Thus, let us assume that a derivation $d$ has the form
$$\infer[(\land)]{P_{\phi}(t,s),\Delta}{
\Prf_S(t,\gn{\neg\phi(\num s)}),\Delta & \phi(s),\Delta}
$$
where $\phi\in\Pi^\LL_{n+1}$ and $\gn{\neg\phi(\num s)}$ denotes the result of substituting the term $s$ into the term $\gn{\neg\phi(\num x)}$. Then, by the induction hypothesis, we
obtain some derivations in $U'$ of the formulas
\begin{equation} \label{lp}
\Prf_S(t,\gn{\neg\phi(\num s)})\lor\bgvee\Delta^-
\end{equation}
and
\begin{equation} \label{rp}
\phi(s)\lor\bgvee\Delta^-.
\end{equation}
Since $\Delta^-$ consists of $\Pi^\LL_{n+1}$-formulas, the reflection rule is
applicable to \refeq{rp} and we obtain a $U'$-proof of
$$\tR_{S,n} (\phi(s(\num{\vec x}))\lor\bgvee\Delta^-(\num{\vec x})).$$
Using Lemma \ref{Rn} we derive:
\begin{enumerate}
\item $\tR_{S,n}(\bgvee\Delta^-(\num{\vec x}))\lor \tR_{S,n}(
\phi(s(\num{\vec x})) \qquad$
\item $\bgvee\Delta^-(\vec x) \lor \tR_{S,n}(\phi(s(\num{\vec x}))) \qquad$
\item $\bgvee\Delta^-(\vec x) \lor \tR_{S,n}(\phi(\num s)) \qquad $
(by Lemma 4.1 of \cite{Bek99b})
\item $\bgvee\Delta^-(\vec x) \lor \Diamond_S\phi(\num s)$.
\end{enumerate}
On the other hand, replacing $t$ by an existential quantifier from \refeq{lp} we obtain
$$\bgvee\Delta^-(\vec x)\lor \Box_S\neg\phi(\num s).$$ Together with
4.~by the rule of cut this yields a $U'$-proof of $\bgvee\Delta^-$.
This concludes the proof of Lemma \ref{genz} and thereby of Theorem \ref{redt2}.
\ep

Now we show that under certain conditions one can characterize the closure of $U$ under the rule $\phi/\tR_{S,n}(\phi)$ for $\phi\in\Pi_{n+1}^\LL$ by $\gw$-iterated reflection principles.

Let a signature $\LL$ now be finite. We begin by stating the well-known fact that the schema $I\Delta_0^\LL$ is finitely axiomatizable over $\EA$. This is a corollary of the following proposition (see Lemma 4.2 in \cite{EnPakh19} for a short proof) that we state for a single unary predicate $P$. 
Recall that $x\in_A z$ means that the $x$-th bit in the binary expansion of $z$ is 1. 

\bpr  \label{idelta-fin} The following are equivalent over $\EA$:
\benr\item $I\Delta_0(P)$;
\item $\al{n}\ex{z}\al{x<n} (P(x)\eqv x\in_A z)$.
\eenr
\epr 

\bcor If $\LL$ is finite, $I\Delta_0^\LL$ is finitely axiomatizable over $\EA$.\ecor 

The following theorem is crucial for several results in this paper. 

\bt \label{ref-fin} Let $\LL$ be finite. Then 
the schema $\tRFN{\Pi_{n+1}^\LL}{S}$ is finitely axiomatizable over $\EA$.
\ignore{\item If $S$ contains $\UTB_\LL$, the schema $\tRFN{\Pi_n^\LL(\TT)}{S}$ is finitely axiomatizable over $\EA+\UTB_\LL$.
\eenr }
\et

\bp\ 
Let $\pi\in\Pi_1^\LL$ denote the formula axiomatizing  $I\Delta_0^\LL$ over $\EA$. We know that the schema $\tR_{S,n}$ implies $\pi$.
We claim that $\tR_{S,n}$ is equivalent over $\EA$ to a conjunction of $\pi$ and the formula 
\beq \label{rf-tr1} \al{\phi\in \Pi_{n+1}^\LL}(\Box_{S}\phi \to \Tr_{n+1}(\phi)),
\eeq 
where $\Tr_m$ is the truth definition for $\Pi_m^\LL$-formulas constructed in Appendix~A. 
Clearly, \refeq{rf-tr1} implies $\tR_{S,n}$ in $\EA+\pi$: For each $\phi\in\Pi_{n+1}^\LL$ we can infer from $\Box_{S}\phi(\num x)$ the formula  $\Tr_{n+1}(\phi(\num x))$ by \refeq{rf-tr1} and then $\phi(x)$ using $\pi$. 

For the opposite implication we first remark that $I\Delta_0^\cL$ and $\pi$ are provable from $\tRFN{\Pi_{n+1}^\LL}{S}$ by Lemma~\ref{ref-ind}. In order to prove \refeq{rf-tr1} we consider two cases: $n=0$ and $n>0$. For $n=0$ we use  Theorem \ref{tr-pi-m} (ii) saying that, for all $\phi\in\Pi_{1}^\LL$, $$\EA\vdash\phi\to \Tr_{1}(\phi).$$
This fact is formalizable in $\EA$. Then, since $S$ provably contains $\EA$, $\Box_S\phi$ implies $\Box_S\Tr_{1}(\phi)$ and hence $\Tr_{1}(\phi)$ by reflection.  

If $n>0$ then we reason in $\EA$ as follows. If $\phi\in\Pi_{n+1}^\cL$ and 
$\Box_S\phi$ then $\Box_S(\pi\to\phi)$. Since $\pi$ provably entails $I\Delta_0^{\cL}$ over $S$, and by Theorem \ref{tr-pi-m} (i) $$\EA^{\cL}\vdash \phi\eqv \Tr_{n+1}(\phi),$$
we obtain $\Box_S(\pi\to\Tr_{n+1}(\phi))$. The formula $\pi\to\Tr_{n+1}(\phi)$ is logically equivalent to a $\Pi_{n+1}^\cL$-formula. Hence, by $\Pi_{n+1}^\cL$-reflection, we can infer $\pi\to\Tr_{n+1}(\phi)$. By an application of $\Pi_{1}^\cL$-reflection we prove $\pi$ and thus we can infer $\Tr_{n+1}(\phi)$.
\ep 

Let $\tR_{S,n}^0:= \top$, $\tR_{S,n}^{k+1}:=\tR_{S,n}({\tR_{S,n}^{k}})$. We can then define $\tR_{S,n}^\gw$ as the (infinite) schema $\{\tR^k_{S,n}:k<\gw\}$.

\bl\ \label{finax}
Suppose $S$ provably contains $U$ and $n\geq 1$. The following theories are equivalent:
\benr
\item $U+\phi/\tR_{S,n}(\phi)$ for $\phi\in\Pi^\LL_n$;
\item $U+\{\tR^k_{S,n}:k<\gw\}$.
\eenr
\el
\bp\ Since the schemata $\tR_{S,n}(\psi)$ are finitely axiomatizable, by external induction on $k$ we can derive $\tR_{S,n}^k$ using $k$ applications of the rule. Hence, theory (i) contains (ii).

The theory $U+\{\tR_{S,n}^k:k<\gw\}$ is closed under the rule: if $U+\tR_{S,n}^k\vdash \phi$, then $S+\tR_{S,n}^k$ contains $S+\phi$, therefore $U+\tR_{S,n}(\tR_{S,n}^k)\vdash \tR_{S,n}(\phi)$.\ep

We now combine Theorem \ref{redt2} and Lemma \ref{finax} into the following theorem that will be used below.

\bt \label{redt3} Suppose $\LL$ is finite. If $U$ is a $\Pi^\LL_{n+2}$-axiomatizable extension of $\EA$ and $S$ provably contains $U$, then
$U+\tR_{S,n+1}$ is a $\Pi^\LL_{n+1}$-conservative extension of the theory $U+\{\tR_{S,n}^k:k<\gw\}$.
\et

\brem The proofs of Theorems \ref{redt2} and \ref{redt3} are formalizable in $\EA^+$, which leads to $\EA^+$-provable conservativity of the respective pairs of theories. The essential ingredient that required the use of superexponentiation axiom was the  application of the cut-elimination theorem for first order logic. Well-known superexponential lower bounds on the speed-up of proofs of $\Pi_2$-statements of $I\Sigma_1$ w.r.t.\ $\PRA$ show that the use of superexponentiation axiom here is really necessary~\cite{Pud85, Ignj}.
\erem 

As another corollary we obtain a new proof of a theorem due to Henryk Kotlarski~\cite{Kot86} characterizing arithmetical consequences of global reflection. In a recent paper, M.~{\L}e{\l}yk~\cite{Lel17} proved that $\Delta_0(\TT)$-induction is equivalent to the global reflection principle $\tRFN{\TT}{\EA}$ over the extension of $\EA$ by the full compositional axioms for truth. (This result was also claimed by Kotlarski, but later a gap was found in his proof of reflection by induction.) A somewhat more general result is as follows.

\begin{corollary}[Kotlarski theorem] Let $U$ be an r.e.~extension of $\EA$ in the language of $\EA$. Then
$U+\CT + I\Delta_0(\TT)$ is conservative over $\EA + \RFN^\gw(U)$ for arithmetical sentences.
\end{corollary}
\bp\ Let $\UTB$ denote $\UTB_\LL$ where $\LL$ is the language of $\EA$, and let $S$ be $U+\UTB$. Clearly, $S$ is provably conservative over $U$. Note that $\UTB + \tRFN{\Sigma_1(\TT)}{S}$ contains $S$ and is sufficient to derive all compositional axioms as well as $I\Delta_0(\TT)$, whence
$$
U+\CT + I\Delta_0(\TT) \subseteq \UTB + \tRFN{\Sigma_1(\TT)}{S} \subseteq_{\Pi_1(\TT)} \UTB + \tR_{S,1}^\omega,
$$
where the latter conservation holds by Theorem \ref{redt3} for the language $\LL(\TT)$. Then one shows that the theory
$\UTB + \tR_{S,1}^n$ is conservative over $\EA + \RFN^n(U)$ for arithmetical sentences
by induction on $n$ using Corollary \ref{cor:rfncons} (formalized in $\EA$).
We obtain that $\UTB + \tR_{S,1}^\omega$ is conservative over $\EA + \RFN^\omega(U)$, whence the result follows.
\ep

\paragraph{Open question.} By the results of Kotlarski and {\L}e{\l}yk the theory $\EA+\CT+I\Delta_0(\TT)$ is equivalent to $$\EA+\CT+\tRFN{\TT}{\EA}\equiv \EA+\CT+\tRFN{\Delta_0(\TT)}{\EA+\UTB}.$$ Is this theory equivalent to $\EA+\UTB+\tRFN{\Sigma_1(\TT)}{\EA+\UTB}$?

\section{Theories of iterated truth and reflection} \label{it-tr}

Consider the language $\LL_\alpha := \LL \cup \{\TT_\beta \mid \beta < \alpha \}$, where $\LL$ is the arithmetic language (or its extension by finitely many predicate symbols). We assume fixed an elementary well-ordering representing ordinals up to $\ga$. This ordering determines a natural G\"odel numbering of all objects of $\LL_\ga$.

We interpret $\TT_\beta$ as the truth definition for the language $\LL_\beta$ (note that $\LL_0$ is the language of arithmetic).
For each $\ga$ we define an $\LL_{\alpha + 1}$-theory $\UTB_\alpha$ as $\UTB_{\LL_\ga}$. Also define
$$
\UTB_{<\alpha} := \bigcup_{\beta < \alpha} \UTB_\beta, \qquad \UTB_{\leqslant \alpha} := \UTB_{<\alpha} + \UTB_\alpha.
$$

Even though the language of $\EA+\UTB_{\leqslant\ga}$ is, in general, infinite, this theory can be considered as a definitional extension of a theory $\EA+\UTB_{\leqslant\ga}^*$ formulated in the language $\LL(\TT_\ga)$ with a single truth predicate.

Indeed, for each formula $\phi\in\LL_{\ga+1}$ let $\phi^*$ denote the result of substitution $\TT_\gb(t)/\TT_\ga(\gn{\TT_\gb(\num t)})$, for all subformulas $\TT_\gb(t)$ of $\phi$ and all $\gb<\ga$. This is obviously a first-order interpretation of $\LL_{\ga+1}$ into $\LL(\TT_\ga)$ preserving $\LL(\TT_\ga)$-formulas.  
\ignore{$\phi[\psi(\vec x)/\TT_\ga(\gn{\psi(\num{\vec x})})]$, for all maximal $\LL_\ga$-subformulas $\psi$ of $\phi$. Since every subformula of $\phi$ of the form $\TT_\gb(t)$ is contained in a maximal one, we have $\phi^*\in\LL(\TT_\ga)$.}
Let $\UTB_{\leqslant\ga}^*$ be axiomatized by $\{\phi^*:\phi\in \UTB_{\leq\ga}\}$. Then, the following lemma is easy to verify.

\bl \label{utbsuc} For all $\phi\in\LL_{\ga+1}$,
\benr
\item $\EA+\UTB_\ga\vdash \phi\eqv \phi^*$;
\item $\EA+\UTB_{\leqslant\ga}\vdash \phi$ iff $\EA+\UTB_{\leqslant\ga}^*\vdash\phi^*$.
\eenr
\el
We remark that $\UTB_{\leqslant\ga}^*$ has a $\Pi_2^\LL(\TT_\ga)$-axiomatization, since $(\cdot )^*$ maps $\Pi_2^{\LL_{\ga+1}}$-formulas to $\Pi_2^\LL(\TT_\ga)$-formulas.

The following lemma is easy to prove by a model-theoretic argument, however we need a proof formalizable in $\EA$. Such a proof is only slightly longer and based on a standard idea.

\bl \label{utb_cons_trn}
$\EA+\UTB_{<\ga}$ conservatively extends $\EA+\UTB_{<\gb}$ for $\gb<\ga$.
\el

\bp\ It is sufficient to construct a non-relativizing interpretation of any finite subtheory of $\UTB_{<\ga}$ in $\UTB_{<\gb}$. Consider such a fragment $U$. Let $\ga_1<\ga_2<\dots <\ga_n$ be all the indices of truth predicates occurring in $U$ with $\ga_i\geq\gb$. It is easy to translate $\TT_{\ga_i}$ into the language $\LL_\gb(\TT_{\ga_1},\dots,\TT_{\ga_{i-1}})$ by case distinction:
$$\TT_{\ga_i}^*(x):=\bigvee_{j<k_i} \ex{\vec y}(x=\gn{\phi_j(\num{\vec y})}\land \phi_j(\vec y)).$$
Here $\phi_j(\vec y)$, for $j<k_i$, are all the formulas for which the Tarski biconditionals for $\TT_{\ga_i}$ occur in $U$. Clearly, the translations of all these biconditionals are provable just from the axioms of $\EA$ (in the extended language). Let us denote this interpretation $K_i$.

Now one can argue by induction on $i$ and prove that the part of $U$ in the language $\LL_\gb(\TT_{\ga_1},\dots,\TT_{\ga_{i-1}})$ is interpretable in $\EA+\UTB_{<\gb}$. This is clear for $i=0$. Assuming that $K$ is such an interpretation for $i$, consider the composition of $K_i$ and $K$.
\ep

Let us now fix an elementary well-ordering $(\gL,<)$ and let $\UTB:=\UTB_{<\gL}$. This ordinal notation system can be extended to a slightly larger segment of ordinals up to $\gw(1+\gL)$, e.g., by encoding ordinals $\gw\ga+n$ as pairs $\la\ga,n\ra$. We introduce the following classes of formulas corresponding to the levels of the hyperarithmetical hierarchy up to $\gw(1+\gL)$:\footnote{According to this definition $\Pi_{1+\ga}$ corresponds to $\Pi_1(\mathbf{0}^{(\ga)})$-sets. Many of the formulas below would be simpler if $\Pi_n$ would denote $\Pi_{1+n}$, but we chose to stick to the standard notation.} 
\bi  
\item $\Pi_n:=\Pi_n^{\LL}$ if $n<\gw$;
\item $
\Pi_{\omega(1+\alpha) + n} := \Pi^{\LL_\ga}_{n+1}(\TT_\alpha);$ 
\item $\Pi_{< \lambda} := \bigcup_{\alpha < \lambda} \Pi_\alpha \text{  for limit } \lambda.$
\ei  

Even though the classes $\Pi_\ga$ are formulated, generally speaking, in an infinite language, we can often restrict the language to a single truth predicate. From Lemma \ref{utbsuc} we obtain
\bl \label{pialpha-fin}
\benr
\item Each $\Delta_0^{\LL_\ga}(\TT_\ga)$-formula is equivalent to a $\Delta_0^\LL(\TT_\ga)$-formula in $\EA+\UTB_{\ga}$;
\item
Each $\Pi_n^{\LL_\ga}(\TT_\ga)$-formula is equivalent to a $\Pi_n^\LL(\TT_\ga)$-formula in $\EA+\UTB_{\ga}$;
\item Each $\LL_\ga$-formula is equivalent to a $\Pi_{<\gw(1+\ga)}$-formula in $\EA+\UTB_{<\ga}$.
\eenr
\el

We define the reflection operators, for all $\ga,\gl<\gw(1+\gL)$, $\gl\in\Lim$, as follows:
\begin{eqnarray*}
\tR_{\alpha}(S) & := & \tRFN{\Pi_{1+\alpha}}{S}, \\
\tR_{<\gl}(S) & := & \text{$\tRFN{\Pi_{<\gl}}{S}$ if $\gl\in\Lim$.}
\end{eqnarray*}

Note that for $n <\gw$ we obtain the usual arithmetical reflection principles
$
\tR_n(S) \equiv \RFN_{\Pi_{n+1}}(S).
$
Further, by Lemma \ref{pialpha-fin} we have

\bpr \label{fin-r} 
\benr
\item If $S$ provably contains $\EA+\UTB_\ga$, then over $\EA+\UTB_\ga$
$$\tR_{\omega(1+\ga)+n}(S) \equiv \tRFN{\Pi^{\LL}_{n+1}(\TT_\ga)}{S};$$
\item If $S$ provably contains $\EA+\UTB_\ga$ and $\gb=\gw(1+\ga)+n$, then $\tR_\gb(S)$ is finitely axiomatizable over $\EA+\UTB_\ga$;
\item If $S$ provably contains $\EA+\UTB_{<\ga}$, then over $\EA+\UTB_{<\ga}$ $$\tR_{<\omega(1+\ga)}(S)\equiv \tRFN{\LL_\ga}{S}\equiv \{\tR_\gb(S) \mid \gb < \gw(1+\ga) \}.$$
\eenr
\epr 
We only remark that Statement (ii) follows from (i) and Theorem \ref{ref-fin}. 

By $\equiv_\alpha$ and $\equiv_{<\lambda}$ we denote conservativity for $\Pi_{1+\alpha}$-sentences and $\Pi_{<\lambda}$-sentences, respectively.
The following conservation results obtained from Theorems \ref{th:llpicons} and \ref{redt3} hold provably in $\EA^+$ and together play the main technical role in our treatment.

\begin{theorem} \label{reduction-lim} Let $\gl=\gw(1+\ga)$ and $S$ provably contain $\EA + \UTB_\alpha$. Over $\EA+\UTB$,
\ignore{$\tR_\lambda(S) \equiv_{<\lambda} \{\tR_\alpha(S) \mid \alpha < \lambda \}$} 
$\tR_\lambda(S) \equiv_{<\lambda}\tR_{<\lambda}(S)$.
\end{theorem}
\bp\ Firstly, $\UTB+\tR_\lambda(S)$ is conservative over $\UTB_{\leq\ga}+\tR_\gl(S)$. Since $\UTB_{<\alpha}$ is axiomatized by $\LL_\alpha$-formulas, Theorem~\ref{th:llpicons} is applicable. Therefore, $\EA+\UTB_{\leqslant \alpha} + \tR_{\gl}(S)$ is $\LL_\alpha$-conservative over $\EA+\UTB_{<\alpha} + \tRFN{\LL_\alpha}{S}$.
\ep

We remark that the same theorem works over any extension of $\UTB$ by $\Pi_{<\gl}$-sentences. 
 
\begin{theorem} \label{reduction-suc} Let $V$ be a $\Pi_{1+\ga+1}$-axiomatized extension of $\EA+\UTB$ and let $S$ contain $V$. Then, over $V$,
$\tR_{\alpha + 1}(S) \equiv_\alpha \{\tR_\alpha(S), \tR_\alpha(S + \tR_\alpha(S)), \dots  \}$.
\end{theorem}
\bp\  Let $U\subseteq \Pi_{1+\ga+1}$ be such that $V=\UTB+U$. If $\ga$ is finite, then $V+\tR_{\alpha + 1}(S)$ is an $\LL$-conservative extension of  $U+\tR_{\alpha + 1}(S)$. Hence, the result follows from Theorem \ref{redt3} for $\LL$, which amounts to the usual reduction property in arithmetic \cite[Theorem 2]{Bek99b}.

Suppose $\ga=\gw(1+\gb)+n$, for some $\gb$ and $n$, then $V+\tR_{\alpha + 1}(S)$ is a conservative extension of  $U+\UTB_{\leq\gb}+\tR_{\alpha + 1}(S)$. Theorem \ref{redt3} does not directly apply in this situation, since the language of $\UTB_{\leq\gb}$ might be infinite. However, by Lemma \ref{utbsuc}, $\UTB_{\leq\gb}$ is a conservative extension of a theory $U':=\UTB^*_{\leq \gb}$ formulated in the language $\LL(\TT_\gb)$. Moreover,
$U'$ has a $\Pi^\LL_{2}(\TT_\gb)$-axiomatization. By Theorem \ref{redt3}, $U+U' + \tR_{\alpha + 1}(S)$ is a $\Pi_{1+\ga}$-conservative extension of
$U+U'+ \{\tR_\alpha(S), \tR_\alpha(S + \tR_\alpha(S)), \dots  \}$. Since $U'$ is contained in $\UTB$, we obtain the result.
\ep

\section{Reflection calculus} \label{rcl}

We refer the reader to a note \cite{Bek18b} for a short introduction to strictly positive logic sufficient for the present paper and to \cite{KKTWZ-arx} for more information from a general  algebraic perspective. For background on modal logic and provability logic we refer to the books~\cite{ChZa,Smo85,Boo93}.

\subsection{The system $\Rcl$}
Fix an ordinal (notation system) $(\gL,<)$. Consider a modal language with propositional variables
$p,q$,\dots , a constant $\top$ and connectives $\land$ and $\ga$,
for each ordinal $\ga<\gL$ (understood as diamond modalities).
Strictly positive formulas (or simply \emph{formulas}) are built up
by the grammar:
$$A::= p \mid \top \mid (A\land A) \mid \ga A, \quad \text{where $\ga<\gL$.}$$
\emph{Sequents} are expressions of the form $A\vdash B$ where $A,B$
are strictly positive formulas. The system $\Rcl$ is given by the
following axioms and rules:

\ben
\item $A\vdash A; \quad A\vdash\top; \quad$ if $A\vdash B$ and $B\vdash C$ then $A\vdash C$;
\item $A\land B\vdash A; \quad A\land B\vdash B; \quad$ if $A\vdash B$ and $A\vdash C$ then
$A\vdash B\land C$;
\item  if $A\vdash B$ then $\ga A\vdash \ga B$;\item $\ga\ga A\vdash \ga A$;
\item $\ga A\vdash \gb A \text{ for $\ga>\gb$};$
\item \label{six} $\ga A\land \gb B\vdash \ga(A\land \gb B)$ for
$\ga >\gb$. \een

The system $\Rc_\gw$ is the familiar system $\Rc$ introduced in an equational logic format by Dashkov \cite{Das12en}, the present formulation is from \cite{Bek12a}.\footnote{The system $\Rc\gw$ considered in \cite{Bek14} is a proper extension of $\Rc$ and is different both from $\Rc_\gw$ and $\Rc_{\gw+1}$.} Dashkov showed that $\Rc$ axiomatizes the set of all sequents $A\vdash B$ such that the implication $A\to B$ is provable in Japaridze's polymodal provability logic GLP. Moreover, unlike GLP itself, $\Rc$ is polytime decidable and enjoys the finite frame property (whereas GLP is Kripke incomplete) \cite{Das12en}.

The system $\Rcl$ is a straightforward generalization of $\Rc$ to transfinitely many modalities. It relates to a version of GLP with transfinitely many modalities $\Glp_\gL$ in the same way as $\Rc$ relates to GLP. The system $\Glp_\gL$ was introduced in \cite{Bek05a} and further studied in several papers by Joost Joosten and David Fern\'andez--Duque~(see \cite{JF13,JF14}).

If $L$ is a strictly positive logic, we write $A\vdash_L B$ for the statement that the
sequent $A\vdash B$ is provable in $L$, and $A=_L B$ stands for $A\vdash_L B$ and $B\vdash_L A$. When context allows we will often omit the subscript.

Notice that $\Rc_\gL$ proves the following \emph{polytransitivity} principles: if
$\ga \geq \gb$ then $\ga\gb A\vdash \gb\gb A\vdash \gb A$ and
$\gb\ga A\vdash \gb\gb A\vdash \gb A$.
Also, the converse of Axiom~\ref{six} is provable in $\Rcl$, so that in fact we have \beq \label{sixeq} \ga(A\land \gb B)=_\Rcl \ga A\land \gb B.\eeq

We also mention the following properties: 
\bi \item If $\gw\leq \Lambda<\Omega$ then $\Rc_\Omega$ conservatively extends $\Rcl$. 
\item Derivability problem in $\Rcl$ is polytime reducible to the problem of comparison of ordinal notations in $(\Lambda,<)$. 
\ei 
The first claim can be proved by an easy syntactic argument, see also \cite{BFJ14} where it is done for $\Glp_\Lambda$. The second claim is based on the result of Dashkov~\cite{Das12en} on the polytime decidability of $\Rc$. Since we are working in this paper with elementary well-orderings $\Lambda$, we can therefore assume that $\Rcl$ is elementary decidable.

\subsection{Variable-free fragment of $\Rcl$}\label{closed_fragment_of_RC}
Let $\Rc^0_\gL$ denote the fragment of $\Rcl$ without propositional variables. Formulas of $\Rc_\gL^0$ will serve for us as canonical ordinal notations. This has been studied quite carefully, so we only briefly recall some basic facts all of which can be found in \cite{Bek05a}.

Let $\Fo^\gL$ denote the set of all variable-free $\Rcl$-formulas, and let $\Fo^\gL_\ga$ denote its restriction to the signature $\{\gb:\ga\leq \gb<\gL\}$, so that $\Fo^\gL=\Fo^\gL_0$. For each $\ga<\gL$ we define binary relations $<^\gL_\ga$ on $\Fo^\gL$ by
$$A<^\gL_\ga B \iffdef B\vdash_\Rcl \ga A.$$
Obviously, $<^\gL_\ga$ is a transitive relation invariantly defined on the equivalence classes w.r.t.\ provable equivalence in $\Rcl$. Since $\Rcl$ is elementarily decidable, so are both $=_\Rcl$ and all of $<^\gL_\ga$.

Since $\Rc_\Omega$ conservatively extends $\Rcl$ for $\gw\leq\gL<\Omega$, the structure $(\Fo^\gL,<_\ga^\gL)$ is isomorphic to an initial substructure of $(\Fo^\Omega,<_\ga^\Omega)$. Therefore, it will be convenient for us in this section to fix a maximal possible $\Lambda$, that is, to think about $\Lambda$ as the class of all ordinals.\footnote{One can also choose the first uncountable ordinal as $\Lambda$.} Notationally we will not carry the subscript $\Lambda$ any longer and will refer to $\Rcl$ and $\Rc_\gL^0$ simply as $\Rc$ and $\Rc^0$, respectively. 
 
An $\Rc$-formula without variables and $\land$ is called a \emph{word} (or a \emph{worm} in some treatments). In fact, any such formula syntactically is a finite sequence of letters $\ga$
(followed by $\top$). If $A,B$ are words then $AB$ will denote $A[\top/B]$, that is, the word corresponding to the concatenation of these sequences. $A\circeq B$ denotes the graphical identity of  formulas (words).

The set of all words will be denoted $\Wo$, and $\Wo_\ga$ will denote its restriction to the signature $\{\gb:\gb\geq \ga\}$.  The following facts are from \cite{Bek05a}:
\bi
\item Every $A\in \Fo_\ga$ is $\Rc$-equivalent to a word in $\Wo_\ga$;
\item $(\Wo_\ga/{=_\Rc},<_\ga)$ is well-ordered. In fact, all these structures are order isomorphic to the class of all ordinals.
\ei

In order to compute the order types of words we recall the standard Veblen hierarchy (see \cite{Poh89}). Given a class
$X \subseteq \On$ let $\en_X$ denote its enumerating function. Let $X'$ denote the class of
fixed points of $\en_X$, that is, $X' = \{\ga \in \On : \en_X(\ga) = \ga\}$. Define by transfinite
induction on $\ga$ the so-called critical classes:
\begin{eqnarray*}
\Cr_0 & = & \{\gw^{1+\ga} : \ga\in\On\};\\
\Cr_{\ga+1} & = & \Cr'_\ga;\\
\Cr_\gl & = & \bigcap_{\ga<\gl}\Cr_\ga, \quad
\text{if $\gl$ is a limit ordinal.}
\end{eqnarray*}
Let $\phi_\ga$ be the enumerating function of $\Cr_\ga$. In particular, $\phi_0(\ga) = \gw^{1+\ga}$ and $\phi_1$ enumerates the fixed points of $\phi_0$, that is, $\phi_1(\ga) = \ge_\ga$.
Our definition of $\Cr_0$ and $\phi_0$ deviates slightly from the standard one, because
we start counting with $\gw$, not with $1$. However, this does not change the
definitions of $\Cr_\ga$ for $\ga > 0$.

It is easy to verify that for all $\ga$ the classes $\Cr_\ga$ are closed and unbounded,
and that the functions $\phi_\ga$ are increasing and continuous.
The least ordinal $\ga$ such that $\ga\in\Cr_\ga$ is the Feferman--Sch\"utte ordinal $\Gamma_0$.
It can also be characterized as the limit of the sequence $\phi_0(0)$, $\phi_{\phi_0(0)}(0)$, \dots, in other words, as the first ordinal closed under the operation $\ga\mapsto \phi_\ga(0)$.

Let $o(A)$ denote the order type of the word $A$ in $(\Wo/{=_\Rc},<_0)$. If $X$ is a class of words, we denote $o(X):=\{o(A): A\in X\}$. We also denote by $\ga\uparrow A$ the result of replacing each letter $\gb$ of $A$ by $\ga+\gb$ and $A^+:= (1\uparrow A)$.

The following statements allow one to compute $o(A)$ in terms of the Veblen $\phi$ function \cite{Bek05a}.

\ben
\item If $A\circeq 0^n$ then $o(A)=n$.
 If $A\circeq A_1^+ 0 A_2^+ 0\cdots 0 A_n^+$, where not all
$A_i$ are empty, then
\[o(A)= \gw^{o(A_n)}+\cdots +\gw^{o(A_1)}.\]
\item $o(\Wo^+_{\gw^\ga})=\Cr_\ga$, where $\Wo^+_\gb:=\Wo_\gb\setminus \{\top\}$.
\item If $A\neq\top$, $\ga>0$ and $\ga=\gw^{\ga_1}+\cdots +\gw^{\ga_n}$ is in Cantor normal form, then $$o(\ga\uparrow A)=\phi_{\ga_1}(\dots\phi_{\ga_n}(-1+o(A))\dots).$$
This formula is based on item 2.
\een

We note that Joosten and Fern\'andez--Duque gave a different and nice way to relate ordinal notations to a family of ordinal functions and the Veblen hierarchy \cite{JF14}. Their formulas based on hyperexponential functions can be used here instead of 1--3. A general advantage of the kind of proof-theoretic analysis that we are doing is that it is modular: we treat words as ordinal notations themselves, we know that the ordering we use is well-founded and naturally computable. So, the treatment of Veblen functions only serves the purpose of relating these notations to other more familiar systems and for our self-control. The reader can use instead the functions and the treatment of notations by Fern\'andez-Duque and Joosten. 

\subsection{Reflection algebras} 

Fix an ordinal notation system $(\gL,<)$ and the corresponding language $\LL_\gL$ as in Section \ref{it-tr}. We interpret reflection calculus $\Rcl$ in the semilattice $\fG_S$ of G\"odelian extensions of $S$ where $S$ provably contains $\EA+\UTB$. In doing that we technically follow the treatments in \cite{Bek14,Bek18} where the reader can look for additional details.

\ignore{
By a G\"odelian theory we mean a theory in the language $\LL_\gL$ whose set of axioms comes equipped with an elementary ($\Delta_0(\exp)$) formula defining its set of G\"odel numbers in the standard model of arithmetic. With every such theory $S$ we associate its provability predicate $\Box_S(x)$.}

Fix some base theory $S$, which we suppose to be closed under the $\Sigma_1$-collection rule. We write $S_1\leq_S S_2$ iff $S_1$ $S$-provably extends $S_2$, that is, $$S\vdash \al{x}(\Box_{S_2}(x)\to \Box_{S_1}(x)).$$ We write $S_1=_S S_2$ iff $S_1\leq_S S_2$ and $S_2\leq_S S_1$. The set of all G\"odelian theories provably extending $S$ modulo $=_S$ is a lattice $\fG_S$ where the meet $S_1\land_S S_2$ corresponds to the (naturally defined) union of theories. Reflection principles $R_\ga$ defined in Section \ref{it-tr} induce, for each $\ga$, monotone and  semi-idempotent operators acting on $\fG_S$, provided $S$ contains $\EA+\UTB$. 

\bd \emph{Reflection algebra of $S$} is the structure of semilattice with operators $(\fG_S;\land_S,\top_S,(R_\ga)_{\ga<\Lambda})$.
\ed 

Reflection algebras yield a natural interpretation of the language of $\Rcl$:
$\Rcl$-formulas are sent to (equivalence classes of) G\"odelian theories in $\fG_S$ in such a way that $\top$ corresponds to $\top_S$, $\land$
corresponds to the union of theories $\land_S$, and $\ga$ corresponds to $\tR_\ga$, for each $\ga<\gL$.

An \emph{arithmetical interpretation in $\fG_S$} is a map $*$ from $\Rcl$-formulas to $\fG_S$ satisfying the following conditions:

\bi
\item $\top^*=S$; \quad $(A\land B)^*=(A^*\land_S B^*)$;
\item $(\ga A)^*= \tR_\ga(A^*)$, for all $\ga<\gL$.
\ei

Using Lemma \ref{Rn} we verify the following basic theorem.
\bt\ \label{sound-rc} For all formulas $A,B$ of $\Rcl$, if $A\vdash_{\Rcl} B$ then $A^*\leq_S B^*$, for all arithmetical interpretations $*$ in $\fG_S$.
\et

\bp\ The proof goes by induction on the length of a derivation in $\Rcl$. We only  check the nontrivial case of Axiom 6 of $\Rcl$, which is based on the finite axiomatizability of the schemata $\tR_\ga$. 

We show,  for all theories $S_1,S_2$ containing $\EA+\UTB$ and all $\gb<\ga$, $$\tR_\ga(S_1)\land \tR_\gb(S_2)\vdash \tR_\ga(S_1\land \tR_\gb(S_2)).$$ Reasoning informally inside $\EA+\UTB$, assume $\phi\in\Pi_\ga$ and $S_1\land \tR_\gb(S_2)\vdash\phi(\num x)$. Using Theorem \ref{fin-r}, let $\psi\in\Pi_\gb$ denote a sentence equivalent to $\tR_\gb(S_2)$, then $S_1\vdash\psi\to\phi(\num x)$.
By $\tR_\ga(S_1)$ we infer $\psi\to\phi(x)$. Since $\EA+\UTB+\tR_\gb(S_2)\vdash\psi$, we conclude $\phi(x)$.    
\ep

\brem It is natural to ask if $\Rcl$ is complete w.r.t.\ the considered arithmetical interpretation. We believe that for theories $S$ containing Peano arithmetic the positive answer can be obtained by the standard methods of proving completeness of Japaridze's provability logic $\Glp$. However, we have not worked out the details.
\erem 

\section{Conservation results for iterated reflection} \label{ord-an}

In this section we prove our main conservation result, a Schmerl-type formula, related to the ordinal analysis of systems of iterated reflection principles. In the next section we will apply it to subsystems of second order arithmetic. We will follow the treatment in \cite{Bek18} which is quite general and deals with iterations of arbitrary computable, monotone, semi-idempotent operators. All the operators $\tR_\ga$ satisfy these conditions, so the results of Section 5 of \cite{Bek18} apply. We summarize what we need in the following theorem.

Consider any ordinal notation system (that is, an elementary strict pre-wellordering) $(D,\prec,0)$ and any G\"odelian theory $U$ in $\fG_S$.

\bpr 
We can specify an elementary formula $\rho(\ga,\gb,x)$ defining a family of G\"odelian theories $\tR_\ga^\gb(U)$ in $\fG_S$, where $\gb\in D$ and $\ga\in\gL$, satisfying the following conditions provably in $\EA$: $\tR_\ga^0(U) \equiv U $ and, if $\gb\succ 0$,
\beq
\tR_\ga^\gb(U) \equiv \textstyle{\bigcup}\{\tR_\ga(\tR_\ga^\gy(U)):\gy\prec \gb\}. \label{itr}
\eeq
Moreover, theories $\tR_\ga^\gb(U)$ are unique modulo provable equivalence in $\EA$.
\epr

We will take as $(D,\prec,0)$ the pre-wellorderings $(\Wo_\ga,<_\ga,\top)$ for various $\ga$. We denote by $o_\ga(A)$ the (G\"odel number of a) word $A$ considered as an element of this notation system. Accordingly, $\tR_\gb^{o_\ga(A)}(S)$ denotes the corresponding iteration of $\tR_\gb$ applied to $S$.

We now consider an elementary well-ordering $(\Lambda,<)$ and denote by $\IB_{<\Lambda}$, or simply $\IB$, the theory $\EA^++\UTB_{<\Lambda}$. 
We fix a G\"odelian extension $S$ of $\IB$ and work in $\fG_S$. Let $A_S^*$ denote the interpretation of a word $A$ in $\fG_S$.

\bt \label{schmr} Suppose $S$ is a $\Pi_{1+\ga+1}$-axiomatized extension of $\IB$. Provably in $\EA^+$, for all $A\in \Wo_\ga$, 
$$S+A_S^*\equiv_\ga S+\tR_\ga^{o_\ga(A)}(S).$$
\et

\bp\ We will use reflexive induction on $A\in\Wo_\ga$ in $\EA^+$. The reflexive induction hypothesis states
$$\al{B<_\ga A} \Box_{\EA^+} (S+B_S^* \equiv_\ga S+\tR_\ga^{o_\ga(B)}(S) ),$$ which implies
$$\al{B<_\ga A} \tR_\ga (B_S^*)=_{\EA^+} \tR_\ga(\tR_\ga^{o_\ga(B)}(S)).$$
We prove $S+A_S^* \equiv_\ga S+\tR_\ga^{o_\ga(A)}(S)$ using the induction hypothesis and reasoning informally over $S$. We are going to use formalized versions of Theorems \ref{reduction-lim} and \ref{reduction-suc}, which are available in $\EA$ and $\EA^+$, respectively. 

First, we prove that $S+A^*_S$ contains $\tR_\ga^{o_\ga(A)}(S)$. If $B<_\ga A$ then $A\vdash_\Rcl \ga B$ and hence $A^*_S \leq_S \tR_\ga(B^*_S)$ by the soundness theorem. By the induction hypothesis $\tR_\ga(B^*_S)$ implies $\tR_\ga(\tR_\ga^{o_\ga(B)}(S))$. Hence, $\al{B<_\ga A}\ A^*_S \vdash \tR_\ga(\tR_\ga^{o_\ga(B)}(S))$, and $A_S^* \vdash \tR_\ga^{o_\ga(A)}(S)$ by \refeq{itr}.

Second, we prove $\Pi_{1+\ga}$-conservation. Assume $\pi\in\Pi_{1+\ga}$ and $S+A_S^* \vdash \pi$. We consider the following cases according to the first letter in $A$.

\textsc{Case 1.} $A\circeq\ga B$. Then $\tR_\ga(B^*_S )\vdash \pi$ and  $\tR_\ga(\tR_\ga^{o_\ga(B)}(S) )\vdash \pi$ by the reflexive induction hypothesis. Since $B<_\ga A$ we obtain $\tR_\ga^{o_\ga(A)}(S)\vdash \pi$.

\textsc{Case 2.} $A\circeq (\gb+1) B$ with $\gb\geq \ga$. 
We define $\Rcl$-formulas $Q^\gb_k(p)$ by: $$Q^\gb_0(p):=p,\quad  Q^\gb_{k+1}(p):=\gb(p\land Q_k^\gb(p)).$$
Then, by Theorem \ref{reduction-suc}, since $S$ is $\Pi_{1+\ga+1}$-axiomatized over $\UTB$,
$$S+\tR_{\gb+1}(B^*_S)\equiv_\gb S+\textstyle\bigcup\{ Q^\gb_k(B)^*_S: k<\gw\}.$$
It follows that $S+Q^\gb_k(B)^*_S\vdash \pi$, for some $k$. We have $Q^\gb_k(B)<_\ga A$, hence $\gy:=o_\ga(Q^\gb_k(B))<o_\ga(A)$. By the reflexive induction hypothesis we obtain
$\tR_\ga(\tR_\ga^\gy(S) )\vdash \tR_\ga(Q_k^\gb(B)_S^*)\vdash \tR_\ga(\pi)\vdash \pi$. It follows that $S+\tR_\ga^{o_\ga(A)}(S) \vdash \tR_\ga(\tR_\ga^\gy(S))\vdash \pi$, q.e.d.

\textsc{Case 3.} $A\circeq \gl B$ where $\gl>\ga$ and $\gl\in\Lim$. It follows that $\gl>\ga+1$ and $S$ is axiomatized by $\Pi_{<\gl}$-sentences over $\UTB$. By Theorem \ref{reduction-lim}, over $S$,
$$\tR_\gl(B^*_S)\equiv_{<\gl} \textstyle\bigcup\{\tR_\gb(B^*_S):\ga\leq\gb<\gl\}.$$
Hence, there is a $\gb$ such that $\ga\leq\gb<\gl$ and $\tR_\gb(B^*_S)\vdash \pi$. It follows that  $\tR_\ga(\tR_\gb(B^*_S))\vdash \tR_\ga(\pi)\vdash \pi$. Since $\gb B<_\ga A$ we obtain $\tR_\ga^{o_\ga(A)}(S) \vdash \tR_\ga(\tR_\ga^{o_\ga(\gb B)}(S))$. By the reflexive induction hypothesis $$\tR_\ga(\tR_\ga^{o_\ga(\gb B)}(S))\vdash \tR_\ga(\tR_\gb(B^*_S))\vdash \pi.$$
\ep

This theorem is the main result of this paper. It allows one to easily calculate proof-theoretic ordinals and conservativity spectra for many theories that can be related to progressions of iterated reflection principles of the considered kind. This is explained in more detail in the next section.  

\section{Proof-theoretic analysis by iterated reflection} \label{anref}

\subsection{Stating ordinal analysis results}\label{anref_1}
Theorem \ref{schmr} provides a way to obtain several types of results usually called \emph{proof-theoretic analyses} of formal systems of predicative strength. The type of results for a given system $T$ one is traditionally interested in include:
\benr
\item Consistency proof for $T$ by transfinite induction along a natural elementary well-ordering;
\item Characterization of the class of provably total computable functions of $T$;
\item Characterization of the order types of elementary (or primitive recursive) well-orderings whose well-foundedness is provable in $T$. 
\eenr

It is well-known that all these questions can be reduced to the one of characterizing the set of consequences of $T$ of an appropriate logical complexity class in terms of progressions of iterated reflection principles. 
More specifically, proof-theoretic analysis of a theory $T$ by iterated reflection can be stated as a conservation result
\beq T\equiv_\ga \tR_\ga^\gb(S), \label{cons-itr}\eeq 
for a suitably weak initial theory $S$ (such as $\EA$ or $\IB$) and appropriate ordinal notations $\ga$ (such as $\ga=0,1$ or $\gw$) and $\gb$. If the relation \refeq{cons-itr} holds, then $T$ is called $\Pi_{1+\ga}$-regular and $\gb$ is called its $\Pi_{1+\ga}$-ordinal (details are spelled out below). Naturally occurring theories $T$ are, indeed, usually regular. 

Theorems of the form \refeq{cons-itr} appeared for the first time in the work of Ulf Schmerl~\cite{Schm,Sch82} who showed among other things that, for any $n<\gw$ and a natural ordinal notation system for $\ge_0$,  $$\PA\equiv_n \tR_n^{\ge_0}(\PRA).$$ The method was further developed in \cite{Bek99b}.
Here we briefly recall basic relationships with the more traditional kinds of proof-theoretic analyses associated with logical complexity levels $\Pi_1^0$, $\Pi_2^0$ and $\Pi_1^1$. 

\medskip
\emph{Consistency proofs}, as in (i), can be obtained from a characterization of arithmetical $\Pi^0_1$-consequences of $T$ in terms of iterated consistency assertions. For example, if $T\equiv_0 \tR_0^\ga(\EA)$ holds provably in $\EA^+$ and $\ga$ is a limit ordinal, then 
$$\EA^+\vdash \Con(T)\eqv \al{\gb<\ga}\tR_0^\gb(\EA).$$
Since $\EA^+$ proves $\RFN_{\Sigma_1}(\EA)$, we have 
$$\EA^+\vdash \al{\gb}(\al{\gy<\gb}\tR_0^\gy(\EA)\to \tR_0^\gb(\EA)).$$
Hence, transfinite induction up to $\ga$ over $\EA^+$ is sufficient to derive $\al{\gb<\ga}\tR_0^\gb(\EA)$ and $\Con(T)$. In fact, this argument only requires the use of transfinite induction rule for a specific $\Pi_1$ (or even a $\Delta_0$) formula. 

\medskip
A characterization of the class $\cF(T)$ of \emph{provably total computable functions of $T$} in terms of the fast growing hierarchy of functions, as in (ii), follows from a characterization of $\Pi^0_2$-consequences of $T$ in terms of iterated uniform $\Pi^0_2$-reflection principles over $\EA$. These claims are elaborated in \cite{Bek99b}, here we state a quick summary. 

With an elementary well-ordering $(\Omega,<)$ (canonical or not) one can associate a hierarchy of fast-growing functions $F_\ga:\nat\to\nat$ for all $\ga<\Omega$:
$$F_\ga(x) := \max\:\{2_x^x+1\}\cup\{F_\gb^{(m)}(n)+1:\gb\prec\ga, \gb ,m,n\leq x\}.$$
We define $\cF_\ga$ as the elementary closure of $\{F_\gb:\gb\prec \ga\}.$ These classes are the same as the standard fast growing hierarchy classes for systems of ordinal notation with a reasonable fundamental sequences assignment (see~\cite{Ros84}), however they are defined even when the latter are not available. Then the following theorem holds~\cite{Bek99b}:

\bpr \label{pi^0_2_oa}\benr \item $\cF(\tR^\ga_1(\EA))=\cF_\ga$;
\item If $S\equiv_1 \tR^\ga_1(\EA)$ then $\cF(S)=\cF_\ga$.
\eenr 
\epr 
 Notice that Claim (ii) follows from (i). 

\medskip 
Now we turn to the characterization of the supremum of order types of $T$-provably well-founded elementary well-orderings, as in (iii), aka \emph{the proof-theoretic $\Pi_1^1$-ordinal of $T$}. 
A recent paper by Pakhomov and Walsh~\cite{PW18} provides a general result characterizing this ordinal in terms of progressions of iterated $\Pi_1^1$-reflection principles in the language of second order arithmetic. Equivalently, it can be stated in terms of iterated full uniform reflection principles in the language of first order arithmetic enriched by a free predicate letter $X$. 

The present setup specifically allows for the  extension of the language of the initial theory by a free predicate letter. Hence, in this context $\Pi_1^1$-analysis amounts to $\Pi_{<\gw}^X$-analysis of a formal system in the language with $X$. The following proposition is a consequence of Theorem~5.9 and Lemma~4.15 of \cite{PW18}.

\begin{proposition} \label{pi^1_1_oa}
\benr \item The $\Pi^1_1$-ordinal of $\tR_{<\gw}^{1+\ga}(\EA^X)$ is $\ge_\alpha$; \item 
If $S\equiv_{<\gw} \tR_{<\gw}^{1+\ga}(\EA^X)$ then the $\Pi^1_1$-ordinal of $S$ is $\ge_\alpha$.
\eenr 
\end{proposition}
As before, Claim (ii) directly follows from Claim (i). 

Summing up our short discussion ordinal analysis we see that it is advantageous \emph{to state} ordinal analysis results in the form \refeq{cons-itr}, for these conservativity relationships yield all the main types of results found in the literature and sometimes, especially in complexity $\Pi_1^0$, provide sharper characterizations. It is a different matter how one can obtain results such as \refeq{cons-itr}.

Conservation results can be obtained by a variety of proof- and model-theoretic methods. However, technical details can sometimes become simpler if one uses methods adapted to deal with reflection principles. The role of reflection algebras is just that, as explained in the next subsection.   

\subsection{The use of reflection algebras for proof-theoretic analysis}
Proof-theoretic analyses such as (i)--(iii) above can be obtained by the following two steps:
\ben
\item  For a given theory $T$ we either find an axiomatization by a finite combination of reflection principles $A^{*}_S$  over a chosen basic system $S$ (such as $\EA$ or $\EA^++\UTB$ in our case) or reduce $T$  to a theory $U$ of this form (that is, prove an  apropriate partial conservation result between $T$ and $U$).
\item  Reduce $A^{*}_S$ to a transfinite iteration of reflection principles of a particular complexity class, such as $\Pi_1$, $\Pi_2$, or $\Pi_{<\gw}^X$. In our case this is achieved by Theorem \ref{schmr} and this ordinal will be equal to $o_\ga(A)$, for some $\ga$, typically $\ga=0$, $1$ or $\gw$, and some $A\in\Wo_\ga$. 
\een 

\ignore{Let an elementary well-ordering be fixed. 
For a theory $T$ whose language contains $\LL_\Lambda$, one can define the notion of $\Pi_{1+\ga}$-ordinal by
$$\ord_\ga(T):= $$}

We remark that Step 1 above is, in practice, relatively easy, since for natural theories $T$ suitable reflection principles are often well-known. The main difficulty in applying these methods is in Step 2, which depends on the overall structure of reflection principles of various complexity levels not to contain big `gaps'. Theorem \ref{schmr} is based on the two main conservation results (Theorems \ref{reduction-lim} and \ref{reduction-suc}) that can be seen as stating that there are no such gaps. 

This idea was first applied to Peano arithmetic and its fragments where Japaridze's polymodal provability logic $\Glp$ was used instead of the reflection calculus $\Rc$ \cite{Bek04,Bek05}. It was later remarked that for these kind of applications the strictly positive fragment of $\Glp$ axiomatized by $\Rc$ is sufficient and allows for a simpler treatment \cite{Bek12a}. 

\subsection{A case study: Analysis of $\ACA$}\label{ACA_case}
This method of analysis, in the simplest situation going beyond Peano arithmetic, can be illustrated by the well-known example of the second order theory $\ACA$. This system extends $\PA$ by the schemata of induction, for all second order formulas, and by the comprehension schema:
\beq\ex{Y}\al{x}(x\in Y\eqv \phi(x)), \label{ca}\eeq for each arithmetical formula $\phi$ (possibly with first- and second-order parameters but not containing $Y$ as a parameter).

It is well-known that $\ACA$ is $\omega$-mutually-interpretable with the extension of $\EA$ in the language $\LL(\TT)$ by the full induction schema together with the compositional axioms for truth (see e.g.~\cite{Halb}). We denote the latter theory $\PA(\TT):=\EA+\CT+I\LL(\TT)$. One way, the interpretation works by defining the domain of set variables by an elementary  formula $\Set(e)$ expressing \emph{``$e$ is the G\"odel number of an arithmetical formula with one free variable''},
and by letting $$x\in e\iffdef \Set(e)\land \TT(e(\num{x})).$$
The other way, the interpretation works by defining in $\ACA$ a $\Delta_1^1$ truth predicate for arithmetical formulas (see, e.g., \cite{Tak} or Appendix \ref{Tr_def_appendix} of this paper). 

By Lemma \ref{ref-ind}, full induction in $\LL(\TT)$ follows from full reflection. Vice versa, the standard argument based on cut-elimination shows that full induction in the presence of $\CT$ proves full reflection. Let $S:=\EA^++\UTB_{\LL(\TT)}$.  We have already noted that the compositional axioms are provable from $\tRFN{\Sigma_1(\TT)}{S}$ over $S$. Thus, $$\PA(\TT)\equiv S+\tRFN{\LL(\TT)}{S}.$$
The latter theory is the same as
$$S+\tR_{<\gw 2}(S)\equiv_{<\gw 2} S+\tR_{\gw 2}(S).$$
Hence, we can apply Theorem \ref{schmr} to $S+\tR_{\gw 2}(S)$. Note that $o_\gw((\gw 2)\top)=o(\gw\top)=\ge_0$. Therefore, Theorem~\ref{schmr} yields that $\PA(\TT)$ is conservative over $S+\mathsf{R}_\gw^{\ge_0}(S)$ for $\Pi_\gw$-sentences. The latter theory is conservative over $\ge_0$ times iterated arithmetical uniform reflection over $\EA^+$ (or equivalently over $\PA$).
This gives a characterization of arithmetical consequences of $\PA(\TT)$ and of $\ACA$. The same reasoning, when carried out in the language with a free predicate letter $X$, yields the well-known bound $\ge_{\ge_0}$ on the provably well-founded elementary orderings in $\ACA$, by \cite{PW18}.

We also obtain that  $o_n((\gw 2)\top)=o((\gw 2)\top))=\phi_1(\phi_1(0))=\ge_{\ge_0}$, for all $n<\gw$. Hence, by Theorem~\ref{schmr} $\PA(\TT)$ is a $\Pi_{n+1}$-conservative extension of  $\EA^++\tR_n^{\ge_{\ge_0}}(\EA^+)$, for each $n<\gw$. This shows, in particular, that $\Pi_1^0$-consequences of $\ACA$ are axiomatized by $\ge_{\ge_0}$ times iterated consistency over $\EA^+$. Similarly, the class of provably total computable functions of $\ACA$ is the $\ge_{\ge_0}$-th class of the extended Grzegorczyk hierarchy.

Similarly, one can analyze subtheories of $\PA(\TT)$ defined by restricted induction. These theories have been studied by Kotlarski and Ratajczyk~\cite{KotRat90} who obtained a characterization of their arithmetical consequences in terms of iterated arithmetical uniform reflection principles. 

By relativization of the usual proof in arithmetic we have that $\EA+\CT+I\Sigma_n(\TT)$ is equivalent to $S+\tRFN{\Pi_{n+2}(\TT)}{S}$. Note that we need $\CT$ to infer reflection from induction, but $\Sigma_1(\TT)$-reflection over $\UTB$ already implies $\CT$. Therefore,
$$\EA+\CT+I\Sigma_n(\TT) \equiv S+\mathsf{R}_{\gw+n+1}(S)\equiv_{\gw} S+\mathsf{R}_{\gw}^{\gw_{n+1}}(S).$$
We have $o_\gw((\gw+n+1)\top)=o((n+1)\top)=\gw_{n+1}$, where $\gw_0=1$, and $\gw_{n+1}=\gw^{\gw_n}$. Hence, $\EA+I\Sigma_n(\TT)$ is conservative over $\gw_{n+1}$-fold iteration of arithmetical uniform reflection over $\EA$.

Now we can remark that, under the standard interpretation of the second order language in $\LL(\TT)$, arithmetical formulas with set parameters are translated into $\Delta_0(\TT)$-formulas, and hence $\Sigma_n^1$-formulas are translated to $\Sigma_n(\TT)$-formulas. Hence, $\EA+\CT+I\Sigma_n(\TT)$ interprets $\ACA_0+I\Sigma^1_n$. Therefore, its $\Pi_{<\gw}$-ordinal is bounded by $\gw_{n+1}$, and $\Pi_{1+m}^0$-ordinal by $o_m((\gw+n+1)\top)=o((\gw+n+1)\top)=\phi_1(\phi_0^{(n+1)}(0))=\ge_{\gw_{n+1}}$.

Applications to theories of iterated arithmetical comprehension will be more comprehensively presented in Section \ref{sec-ord}. 

\subsection{Conservativity spectra}

The notions of $\Pi_{1+n}^0$-ordinal~\cite{Bek99b} and conservativity spectrum~\cite{Joo15a,Bek18b} can be naturally generalized to the classes of the hyperarithmetical hierarchy.

Let us fix, as before, an elementary well-ordering $\Lambda$ and the corresponding language $\LL_\Lambda$. By $\IB_\ga$ we denote the extension of $\EA^+$ by the axioms $\UTB_{\leq\ga}$, and $\IB_{<\ga}$ denotes the union of theories $\IB_\gb$ for all $\gb<\ga$. We will consider iterations of reflection principles along another elementary well-ordering that need not be related to $\Lambda$. Let us fix such a well-ordering $\Omega$. 

\ignore{
and consider another elementary well-ordering $\Omega$.   
Let $S$ be a G\"odelian extension of $\EA^++\UTB_\Lambda$ and a fixed  elementary recursive well-ordering. In this section we additionally assume that $\Omega$ is an epsilon number and is equipped with elementary terms representing the ordinal constants and functions $0,1,+,\cdot,\gw^x$. These functions should provably in $\EA$ satisfy some minimal natural axioms NWO listed in \cite{Bek95}. We call such well-orderings \emph{nice}. Recall the following definitions from \cite{Bek99b} (writing $1$ for $1_{\EA^+}$):
}

\bd Let $S$ be a G\"odelian extension of $\IB$.
\bi
\item \emph{(lower) $\Pi_{1+\ga}$-ordinal of $S$}, denoted $\ord_\ga(S)$, is the supremum of all $\gb\in\Omega$ such that $S\vdash \tR_\ga^\gb(\IB)$;
\item \emph{upper $\Pi_{1+\ga}$-ordinal of $S$}, denoted $\ord^u_\ga(S)$, is the infinum of all $\gb\in\Omega$ such that $S$ is $\Pi_{1+\ga}$-conservative over $\tR_\ga^\gb(\IB)$;
\item $S$ is \emph{$\Pi_{1+\ga}$-regular} if lower and upper  $\Pi_{1+\ga}$-ordinals of $S$ coincide, that is,  $\ord_\ga(S)=\ord^u_\ga(S)$.
\ei
\ed 

As for finite $\ga$, $\Pi_{1+\ga}$-ordinals are insensitive to  $\Pi_{1+\ga}$-conservative extensions and to extensions by consistent $\Sigma_{1+\ga}$-axioms. 

\bpr \label{insense} For for all $\ga<\gw(1+\Lambda)$,
\benr\item If $T$ is $\Pi_{1+\ga}$-conservative over $S$ then $\ord_\ga(S)\geq \ord_\ga(T)$;
\item Suppose $S$ is $\Pi_{1+\ga}$-regular. If $T$ is axiomatized by $\Sigma_{1+\ga}$-sentences and $S\cup T$ is consistent, then $\ord_\ga(S\cup T)=\ord_\ga(S)$.  
    \eenr
\epr
\bp\ We prove the first claim using the fact that $\tR_\ga^\gb(\IB)$ is a $\Pi_{1+\ga}$-axiomatized extension of $\IB$. The second claim follows from the well-known result that $\tR_\ga(U)$ is not contained in any consistent $\Sigma_{1+\ga}$-axiomatized extension of $U$. \ep

The sequence of $\Pi_{1+n}$-ordinals of a given system $S$ is sometimes called its conservativity spectrum. This sequence bears the more detailed information about the strength of a given theory at various levels of logical complexity than any individual proof-theoretic ordinal. Joost Joosten~\cite{Joo15a} studied such sequences for finite $n$. He showed for theories between $\EA^+$ and $\PA$ that their conservativity spectra correspond to decreasing sequences of ordinals below $\ge_0$ of a certain kind, that is, to the points of the so-called Ignatiev frame. Beklemishev~\cite{Bek18} showed that the set of spectra naturally bears the structure of an $\Rc$-algebra with conservativity operators.   
An immediate generalization of conservativity spectra in arithmetic is as follows.

\bd \emph{Conservativity spectrum of $S$} is the sequence $(\ga_\xi)_{\xi<\gw(1+\gL)}$ such that $\ga_\xi=\ord_\xi(S)$, for all $\xi$.
\ed
More explicitly, we can call it the $\gw(1+\gL)$-conservativity spectrum of $S$.

As an example consider the theory $\ACA$ or $\PA(\TT)$. Since $\PA(\TT)$ is formulated in $\cL_{\gw 2}$, we consider its $\gw 2$-spectrum. As we have seen, this is the sequence
$$(\ge_{\ge_0},\ge_{\ge_0},\dots;\ge_0,\ge_0,\dots).$$
Observe that if, for example, instead of $\PA(\TT)$ we consider the set of its arithmetical consequences (axiomatized by $\ge_0$-fold iteration of arithmetical uniform reflection principle), then its conservativity spectrum looks as follows: 
$$(\ge_{\ge_0},\ge_{\ge_0},\dots;0,0,\dots).$$

We conjecture that in general conservativity spectra consist precisely of $\ell$-sequences in the sense of \cite{JF13}. Therefore, they one-to-one correspond the points of the generalized Ignatiev frames constructed by Fern\'andez and Joosten. However, proving this conjecture remains out of the scope of the present paper.

\section{Analysis of second order systems} \label{sec-ord}

In this section we show how Theorem \ref{schmr} can be used to obtain ordinal analysis of some systems of second order arithmetic of `predicative' strength. 

\subsection{Ordinal analysis of iterated arithmetical comprehension} 
For the purposes of this section we use as the ordinal notion system $\Lambda$ any natural elementary ordinal notation system that contains the symbol for $0$, is closed under addition, function $\omega^x$, and the binary Veblen function $\varphi_x(y)$. We also assume that $\mathsf{EA}^+$ is capable to prove standard universal properties of this functions in the ordinal notation system, e.g., that addition is associative, that $\omega^{\alpha}<\omega^{\beta}$ iff $\alpha<\beta$, that $\varphi_{\alpha}(\varphi_{\beta}(\gamma))=\varphi_{\beta}(\gamma)$, for $\beta>\alpha$, etc. Note that it is well-known to be the case for the standard ordinal notation systems for the ordinal $\Gamma_0$.

The base theory of second-order arithmetic we consider is the well-known theory $\mathsf{ACA}_0$, that is,  the extension of $\mathsf{EA}$ by the scheme of arithmetic comprehension \refeq{ca}
and the axiom of set-induction $$0\in X\land \forall x\;(x\in X\to S(x)\in X)\to \forall x\;(x\in X).$$
Due to our choice of base theory, we freely use pseudo-terms $\{x\mid F(x)\}$, when $F(x)$ is a $\Pi^0_\infty$-formula.

The are various ways how one could axiomatized the theory of iterated $\Pi_1^0$-comprehension $(\Pi_1^0\tCA_0)_{\alpha}$.  Officially we will use the axiomatization that extends $\mathsf{ACA}_0$ by the following scheme:
\begin{equation}\label{itcmp_schm}\exists X\forall \beta\le\alpha \;((X)_\beta= \{x\mid F(x,(X)_{<\beta})\}),\mbox{ for $\Pi^0_{\infty}$ formulas $F(x,X)$};\end{equation}
here $(Y)_\beta:=\{x:\la \beta, x\ra\in Y\}$ and $(Y)_{<\beta}:=\{x:\ex{\gamma<\beta} \la \gamma,x\ra\in Y\}$.
However, in certain situations it will be useful for us to use a different axiomatization that extends $\mathsf{ACA}_0$ by the axiom stating that for any set $X$ the $\alpha$-th Turing jump $X^{(\alpha)}$ exists.

We also consider the following theories
\begin{enumerate}
    \item $(\Pi_1^0\tCA)_{\alpha}:=\mathsf{ACA}+(\Pi^0_1\mbox{-}\mathsf{CA}_0)_{\alpha}$;
    \item $(\Pi_1^0\tCA_0)_{<\alpha}:=\mathsf{ACA_0}+\bigcup\limits_{\beta<\alpha}(\Pi^0_1\mbox{-}\mathsf{CA}_0)_{\beta}$;
    \item $(\Pi_1^0\tCA)_{<\alpha}:=\mathsf{ACA}+(\Pi_1^0\tCA_0)_{<\alpha}$.
\end{enumerate}
Note that some usual systems are of this form: $\mathsf{ACA}_0\equiv(\Pi^0_1\mbox{-}\mathsf{CA}_0)_{1}$, $\mathsf{ACA}\equiv(\Pi^0_1\mbox{-}\mathsf{CA})_{1}$, $\mathsf{ACA}_0^+\equiv(\Pi^0_1\mbox{-}\mathsf{CA}_0)_{\omega}$, and $\mathsf{ACA}^+=(\Pi^0_1\mbox{-}\mathsf{CA})_{\omega}$.

Recall that the definitions of the theories $\mathsf{UTB}_{<\alpha}$ and $\mathsf{UTB}_{\le\alpha}$ required to first fix the ground language $\mathcal{L}$ (that should be an extension of arithmetical language by finitely many predicate symbols). In this section it will be sufficient to only consider the case of $\mathcal{L}$ expanding the arithmetical language by only one unary predicate symbol $X(x)$. The inclusion of the free predicate symbol is important for studying the $\Pi^1_1$-consequences of second-order theories. We identify the unary predicate letter $X(x)$ with the set variable $X$ and the atomic formulas $X(t)$ with the atomic formulas $t\in X$. This allows us to identify a $\Pi^1_1$-sentence $\forall Y\varphi(Y)$ with the $\LL$-sentence $\varphi(X)$ and to talk about the sets of $\Pi^1_1$-consequences of theories whose language contains $\LL$. 

We define a translation $\mathcal{E}_X(\cdot)$ of the language $\mathcal{L}_{\Lambda}$ into the language of second-order arithmetic. The set variable $X$ is the only parameter of the translation. The translation $\mathcal{E}_X(\cdot)$ interprets symbols of  arithmetical language by themselves and the predicate $X$ by the set $X$. 

To define the interpretations $\mathcal{E}_X(\TT_{\alpha}(x))$ of truth predicates we will use partial truth predicates encoded by sets. Let $\mathsf{v}(x)$ be the term evaluation function, i.e. the function that maps a G{\"o}del number of any closed term $t$ in the language of arithmetic to its numerical value. A set $P$ of G{\"o}del numbers of $\LL_{\alpha}$ sentences is called a \emph{compositional truth definiton for $\LL_{\alpha}$ over a set $X$} if the following conditions hold:
\begin{enumerate}
    \item $\varphi\in P \mathrel{\leftrightarrow} \varphi$, for any atomic arithmetic sentence $\varphi$;
    \item $X(t)\in P\mathrel{\leftrightarrow} \mathsf{v}(t)\in X$, for for any closed term $t$;
    \item $\TT_{\beta}(t)\in P \mathrel{\leftrightarrow} \mathsf{v}(t)\in \LL_{\beta} \land \mathsf{v}(t)\in P$, for any closed term $t$ and $\beta<\alpha$;
    \item $(\varphi\land \psi)\in P\mathrel{\leftrightarrow} \varphi\in P\land \psi\in P$, for any sentences $\varphi, \psi \in\LL_{\alpha}$;
    \item $(\lnot\varphi)\in P\mathrel{\leftrightarrow} \varphi\not\in P$, for any sentence $\varphi \in\LL_{\alpha}$;
    \item $(\forall x\;\varphi(x))\in P \mathrel{\leftrightarrow} \forall x\; (\varphi(\underline{x})\in P)$, for any $\varphi(x) \in\LL_{\alpha}$;
\end{enumerate}
In this case we write $\mathsf{TP}^X_{\alpha}(P)$.  We finish the definition of  $\mathcal{E}_X$ by interpreting the  predicates $\TT_{\alpha}(x)$ as
$$\exists \beta\le \alpha\;(x\in \LL_{\beta}\land \forall P(\mathsf{TP}^X_{\beta}(P)\to x\in P)).$$

By this translation we will consider $\mathcal{L}_{\Lambda}$ to be a sub-language of the language of second-order arithmetic and will freely use expressions like $S \equiv_{\alpha} U$, where $S$ is a second-order theory and $U$ is an  $\mathcal{L}_{\Lambda}$ theory (or vice versa). In  terms of the translation $\mathcal{E}_X(\cdot)$ this means that for each $\Pi_{1+\alpha}$ sentence $\varphi$ we have $S\vdash \forall X\; \mathcal{E}_X(\varphi) \iff U\vdash \varphi.$

Now let us state our two central results connecting systems of second-order arithmetic with the systems of iterated truth definitions.
\begin{theorem} \label{ica_ref_0} For $\alpha\geq 1$, $$(\Pi^0_1\mbox{-}\mathsf{CA}_0)_{\omega^{\alpha}}\equiv_{<\omega^{\alpha+1}}\mathsf{R}_{<\omega^{\alpha+1}}(\mathsf{IB})$$\end{theorem}
\begin{theorem} \label{ica_ref_ind} For $\alpha\geq 1$, $$(\Pi^0_1\mbox{-}\mathsf{CA})_{\omega^{\alpha}}\equiv_{<\omega^{\alpha+1}+\omega}\mathsf{R}_{<\omega^{\alpha+1}+\omega}(\mathsf{IB}).$$\end{theorem}

Before proving these two theorems we will elaborate on how the proof-theoretic analysis for theories of iterated $\Pi^0_1$-comprehension follows from the two results and the discussion in Section \ref{anref_1}.

Theorems \ref{ica_ref_0} and \ref{ica_ref_ind} cover only the case of theories of $\beta$-times iterated comprehension, where $\beta$ is of the form $\omega^{\alpha}$. But it is easy to reduce more general case to this one. For an ordinal $\alpha$ with the Cantor normal form $\omega^{\alpha_1}+\ldots+\omega^{\alpha_n}$ the theories $(\Pi^0_1\mbox{-}\mathsf{CA}_0)_{\alpha}$ and $(\Pi^0_1\mbox{-}\mathsf{CA})_{\alpha}$ coincide with the theories $(\Pi^0_1\mbox{-}\mathsf{CA}_0)_{\omega^{\alpha_1}}$ and $(\Pi^0_1\mbox{-}\mathsf{CA})_{\omega^{\alpha_1}}$, respectively. 

For the case of limit theories $(\Pi^0_1\mbox{-}\mathsf{CA}_0)_{<\alpha}$ and $(\Pi^0_1\mbox{-}\mathsf{CA})_{<\alpha}$, where $\alpha\ge 1$, by a combination of compactness, Theorem \ref{ica_ref_0}, and Theorem \ref{ica_ref_ind} we have
\begin{equation} \label{lim_th_0} (\Pi^0_1\mbox{-}\mathsf{CA}_0)_{<\alpha}\equiv_{<\beta}\mathsf{R}_{<\beta}(\mathsf{IB})\mbox{, where $\beta=\sup\{\delta\omega\mid \delta<\alpha\}$};\end{equation} \begin{equation} \label{lim_th_ind}(\Pi^0_1\mbox{-}\mathsf{CA})_{<\alpha}\equiv_{<\beta}\mathsf{R}_{<\beta}(\mathsf{IB})\mbox{, where $\beta=\sup\{\delta\omega+\omega\mid \delta<\alpha\}$}.\end{equation}

Combining the observations above with Theorem \ref{schmr} and computations from Section \ref{closed_fragment_of_RC} we get the following Schmerl-style formulas for the theories of iterated $\Pi^0_1$-comprehension:
\begin{theorem} \label{schmr_ic}For $\alpha\geq 1$,
\begin{enumerate}
    \item \label{schmr_ic_1} $(\Pi^0_1\mbox{-}\mathsf{CA}_0)_{\omega^{\alpha}}\equiv (\Pi^0_1\mbox{-}\mathsf{CA}_0)_{<\omega^{\alpha+1}}\equiv_{\beta} \mathsf{R}_{\beta}^{\varphi_{\alpha+1}(0)}(\mathsf{IB})$, for $\beta<\omega^{\alpha+1}$;
    \item \label{schmr_ic_2} $(\Pi^0_1\mbox{-}\mathsf{CA})_{\omega^{\alpha}}\equiv (\Pi^0_1\mbox{-}\mathsf{CA})_{<\omega^{\alpha+1}}\equiv_{\beta} \mathsf{R}_{\beta}^{\varphi_{\alpha+1}(\varepsilon_0)}(\mathsf{IB})$, for $\beta<\omega^{\alpha+1}$;
    \item \label{schmr_ic_3} $(\Pi^0_1\mbox{-}\mathsf{CA})_{\omega^{\alpha}}\equiv (\Pi^0_1\mbox{-}\mathsf{CA})_{<\omega^{\alpha+1}}\equiv_{\beta} \mathsf{R}_{\beta}^{\varepsilon_0}(\mathsf{IB})$, for $\omega^{\alpha+1}\le\beta<\omega^{\alpha+1}+\omega$;
    \item \label{schmr_ic_4} $(\Pi^0_1\mbox{-}\mathsf{CA}_0)_{<\omega^{\lambda}}\equiv_{\beta} (\Pi^0_1\mbox{-}\mathsf{CA})_{<\omega^{\lambda}}\equiv_{\beta} \mathsf{R}_{\beta}^{\varphi_{\lambda}(0)}(\mathsf{IB})$, for limit $\lambda$ and $\beta<\omega^{\lambda}$.    
\end{enumerate}
\end{theorem}
\bp\  Let us first prove Claim~\ref{schmr_ic_1}. By Theorem \ref{ica_ref_0} and Theorem \ref{schmr} we have
$$(\Pi^0_1\mbox{-}\mathsf{CA}_0)_{\omega^{\alpha}}\equiv (\Pi^0_1\mbox{-}\mathsf{CA}_0)_{<\omega^{\alpha+1}} \equiv_{<\omega^{\alpha+1}} \mathsf{R}_{<\omega^{\alpha+1}}(\mathsf{IB})\equiv_{\beta} \mathsf{R}_{\beta}^{o_{\beta}(\omega^{\alpha+1})}(\mathsf{IB}).$$
To finish the proof of Claim~\ref{schmr_ic_1} we calculate
$$o_{\beta}(\omega^{\alpha+1})=o(\omega^{\alpha+1})=\varphi_{\alpha+1}(0),$$
where we use results of Section \ref{closed_fragment_of_RC} and the fact that $\omega^{\alpha+1}=\beta+\omega^{\alpha+1}$.

For Claims~\ref{schmr_ic_2} and \ref{schmr_ic_3} by Theorem \ref{ica_ref_ind} and Theorem \ref{schmr} we have
$$(\Pi^0_1\mbox{-}\mathsf{CA})_{\omega^{\alpha}}\equiv (\Pi^0_1\mbox{-}\mathsf{CA})_{<\omega^{\alpha+1}}\equiv_{<\omega^{\alpha+1}+\omega} \mathsf{R}_{<\omega^{\alpha+1}+\omega}(\mathsf{IB})\equiv_{\beta} \mathsf{R}_{\beta}^{o_{\beta}(\omega^{\alpha+1}+\omega)}(\mathsf{IB}).$$
To prove Claim~\ref{schmr_ic_2} we observe  $\omega^{\alpha+1}+\omega=\beta+\omega^{\alpha+1}+\omega$ whence
$$o_{\beta}(\omega^{\alpha+1}+\omega)=o(\omega^{\alpha+1}+\omega)=\varphi_{\alpha+1}(\varphi_1(0))=\varphi_{\alpha+1}(\varepsilon_0).$$
To prove Claim~\ref{schmr_ic_3} we notice  $\omega^{\alpha+1}+\omega=\beta+\omega$ and thus $o_{\beta}(\omega^{\alpha+1}+\omega)=o(\omega)=\varepsilon_0$.

Let us prove Claim~\ref{schmr_ic_4}. Since $\omega^{\lambda}=\sup\{\delta\omega\mid \delta<\omega^{\lambda}\}=\sup\{\delta\omega+\omega\mid \delta<\omega^{\lambda}\}$, by (\ref{lim_th_0}) and (\ref{lim_th_ind}) we obtain 
$$(\Pi^0_1\mbox{-}\mathsf{CA}_0)_{<\omega^{\lambda}}\equiv_{\beta} (\Pi^0_1\mbox{-}\mathsf{CA})_{<\omega^{\lambda}}\equiv_{<\omega^{\lambda}}\mathsf{R}_{<\omega^{\lambda}}(\mathsf{IB})\equiv_{\beta}\mathsf{R}_{\beta}^{o_{\beta}(\omega^{\lambda})}(\mathsf{IB}).$$
Clearly, $o_{\beta}(\omega^{\lambda})=\varphi_{\lambda}(0)$.\ep

Combining Theorem \ref{schmr_ic} with the results of Section \ref{anref_1} we obtain various other ordinal analysis results.
\begin{corollary}For $\alpha\geq 1$,
\begin{enumerate}
\item Theory $\mathsf{EA}^+$ plus transfinite induction for $\Delta_0$-formulas up to $\varphi_{\alpha+1}(0)$ proves the consistency of $(\Pi^0_1\mbox{-}\mathsf{CA}_0)_{\omega^{\alpha}}$.
\item Theory  $\mathsf{EA}^+$ plus transfinite induction for $\Delta_0$-formulas up to $\varphi_{\alpha+1}(\varepsilon_0)$ proves the consistency of $(\Pi^0_1\mbox{-}\mathsf{CA})_{\omega^{\alpha}}$.
\item For limit $\lambda$ the theory $\mathsf{EA}^+$ plus transfinite induction for $\Delta_0$-formulas up to $\varphi_{\lambda}(0)$ proves the consistency of $(\Pi^0_1\mbox{-}\mathsf{CA}_0)_{<\omega^{\lambda}}$ and $(\Pi^0_1\mbox{-}\mathsf{CA})_{<\omega^{\lambda}}$.
\end{enumerate}
\end{corollary}

\begin{corollary}For $\alpha\geq 1$,
\begin{enumerate}
\item $\mathcal{F}((\Pi^0_1\mbox{-}\mathsf{CA}_0)_{\omega^{\alpha}})=\mathcal{F}_{\varphi_{\alpha+1}(0)}$;
\item $\mathcal{F}((\Pi^0_1\mbox{-}\mathsf{CA})_{\omega^{\alpha}})=\mathcal{F}_{\varphi_{\alpha+1}(\varepsilon_0)}$;
\item $\mathcal{F}((\Pi^0_1\mbox{-}\mathsf{CA}_0)_{<\omega^\lambda})=\mathcal{F}((\Pi^0_1\mbox{-}\mathsf{CA})_{<\omega^\lambda})=\mathcal{F}_{\varphi_{\lambda}(0)}$, for limit $\lambda$.
\end{enumerate}
\end{corollary}

\begin{corollary}For $\alpha\geq 1$,
\begin{enumerate}
\item The $\Pi^1_1$-ordinal of $(\Pi^0_1\mbox{-}\mathsf{CA}_0)_{\omega^{\alpha}}$ is $\varphi_{\alpha+1}(0)$;
\item The $\Pi^1_1$-ordinal of $(\Pi^0_1\mbox{-}\mathsf{CA})_{\omega^{\alpha}}$ is $\varphi_{\alpha+1}(\varepsilon_0)$;
\item For limit $\lambda$, the $\Pi^1_1$-ordinal of  $(\Pi^0_1\mbox{-}\mathsf{CA}_0)_{<\omega^\lambda}$ and $(\Pi^0_1\mbox{-}\mathsf{CA})_{<\omega^\lambda}$ is $\varphi_{\lambda}(0)$.
\end{enumerate}
\end{corollary}

\brem In this section we have not covered the theories $\mathsf{ACA}_0=(\Pi^0_1\mbox{-}\mathsf{CA}_0)_1$ and $\mathsf{ACA}=(\Pi^0_1\mbox{-}\mathsf{CA})_1$. This is due to the fact that we wanted to simplify our proofs and the cases of $\mathsf{ACA}_0$ and $\mathsf{ACA}$ require a separate consideration (we already considered the case of $\mathsf{ACA}$ in Section \ref{ACA_case}). The main reason for this complication is that the theories of iterated truth definitions treat successor stages of hyperarthmetical hierarchy $\Pi_{\alpha+1}$ and limit stages $\Pi_{\lambda}$ differently. 
\erem

\subsection{Constructing a universe of sets from  truth definitions}

In this and next section we prove Theorems \ref{ica_ref_0} and \ref{ica_ref_ind}. The result of this section are the following two lemmas.
\begin{lemma} \label{ica_ref_0_w} For any $\alpha$, $$(\Pi^0_1\mbox{-}\mathsf{CA}_0)_{\omega^{1+\alpha}}\subseteq_{<\omega}\mathsf{R}_{<\omega^{1+\alpha+1}}(\mathsf{IB})$$\end{lemma}
\begin{lemma} \label{ica_ref_ind_w}For any $\alpha$, $$(\Pi^0_1\mbox{-}\mathsf{CA})_{\omega^{1+\alpha}}\subseteq_{<\omega}\mathsf{R}_{<\omega^{1+\alpha+1}+\omega}(\mathsf{IB}).$$\end{lemma}
Note that Lemmas \ref{ica_ref_0_w} and \ref{ica_ref_ind_w} essentially are the parts of Theorems \ref{ica_ref_0} and \ref{ica_ref_ind} that are sufficient to prove the upper bounds for  $\Pi^1_1$, $\Pi^0_2$, and $\Pi^0_1$ ordinals of  $(\Pi^0_1\mbox{-}\mathsf{CA}_0)_{\omega^{\alpha}}$ and $(\Pi^0_1\mbox{-}\mathsf{CA})_{\omega^{\alpha}}$.  


As we will see, there is an essential difference between the conservation results of Lemma \ref{ica_ref_0_w} and  \ref{ica_ref_ind_w}. We prove Lemma \ref{ica_ref_0_w} by a model-theoretic construction. Given a model $\mathfrak{M}\models \mathsf{R}_{<\omega^{1+\alpha+1}}(\IB)$ we construct a model $\mathbf{DEF}_{\omega^{\alpha+1}}(\mathfrak{M})\models (\Pi^0_1\mbox{-}\mathsf{CA}_0)_{\omega^{\alpha}}$, where $\mathfrak{M}$ and $\mathbf{DEF}_{\omega^{\alpha+1}}(\mathfrak{M})$ have the same first-order arithmetical part. This does not lead to an interpretation of $(\Pi^0_1\mbox{-}\mathsf{CA}_0)_{\omega^{\alpha}}$ in $\mathsf{R}_{<\omega^{1+\alpha+1}}(\IB)$. Moreover, such an interpretation is impossible, since $(\Pi^0_1\mbox{-}\mathsf{CA}_0)_{\omega^{\alpha}}$ is finitely axiomatizable and at the same time proves the consistency of any finite subtheory of $\mathsf{R}_{<\omega^{1+\alpha+1}}(\mathsf{IB})$. On the other hand, we prove Lemma \ref{ica_ref_ind_w} by constructing an $\omega$-interpretation $\mathcal{D}_{\omega^{\alpha+1}}$ of $(\Pi^0_1\mbox{-}\mathsf{CA})_{\omega^{\alpha}}$ in $\mathsf{R}_{<\omega^{1+\alpha+1}+\omega}(\mathsf{IB})$. This $\omega$-interpretation is provided by a formalized version of the proof of Lemma \ref{ica_ref_0_w}. 

For a model $\mathfrak{M}$ of $\mathsf{UTB}_{<\alpha}$ we define a model of second-order arithmetic $\mathbf{DEF}_{\alpha}(\mathfrak{M})$. The arithmetical part of $\mathbf{DEF}_{\alpha}(\mathfrak{M})$ coincides with the arithmetical part of $\mathfrak{M}$. The second-order part of $\mathbf{DEF}_{\alpha}(\mathfrak{M})$ consists of all the subsets of $\mathfrak{M}$ definable by $\mathcal{L}_{\alpha}$-formulas with parameters from $\mathfrak{M}$.

\proof{ of Lemma \ref{ica_ref_0_w}.}\; It is enough to show that for any $\mathfrak{M}\models \mathsf{R}_{<\omega^{1+\alpha+1}}(\mathsf{IB})$ the model $\mathbf{DEF}_{\omega^{\alpha+1}}(\mathfrak{M})$ satisfies $(\Pi^0_1\mbox{-}\mathsf{CA}_0)_{\omega^{1+\alpha}}$.  

For an $\mathcal{L}_{\omega^{\alpha+1}}$-formula $\varphi(x)$ with parameters from $\mathfrak{M}$ we  denote by $D_{\varphi}$ the set $\{a\mid \mathfrak{M}\models \varphi(a)\}$. Clearly, each set $A$ from $\mathbf{DEF}_{\omega^{\alpha+1}}(\mathfrak{M})$ is $D_{\varphi}$ for some $\varphi$. Naturally, we can treat $\Pi^0_\infty$ formulas $\psi(\vec{x})$ with parameters from $\mathbf{DEF}_{\omega^{\alpha+1}}(\mathfrak{M})$  as $\mathcal{L}_{\omega^{\alpha+1}}$ formulas: We transform $\psi$ to an $\mathcal{L}_{\omega^{\alpha+1}}$-formula by replacing the subformulas $t\in D_{\varphi}$ with $\varphi(t)$. Moreover, we apply the same treatment to arbitrary $\mathcal{L}_{\Lambda}$-formulas and use the formulas $t\in D_{\varphi}$ as shorthands for $\varphi(t)$. Since we talk about definable hierarchies, in addition we will use formulas $t\in (D_{\varphi})_{\beta}$ as shorthands for $\varphi(\langle t,\beta\rangle)$ and formulas $t\in (D_{\varphi})_{<\beta}$ as shorthands for $\varphi(t)\land \exists \gamma<\beta\exists x\;(t=\langle \gamma,x\rangle)$. 

Since $\mathsf{R}_{<\omega^{1+\alpha+1}}(\mathsf{IB})$ proves induction for $\mathcal{L}_{\alpha}$ formulas, we see that the model $\mathbf{DEF}_{\omega^{\alpha+1}}(\mathfrak{M})$ satisfies  the second-order induction axiom. Clearly, the result of comprehension for a $\Pi^0_\infty$ formula $\varphi(x)$ with parameters from $\mathbf{DEF}_{\omega^{\alpha+1}}(\mathfrak{M})$ is the set $D_{\varphi}$. Thus $\mathbf{DEF}_{\omega^{\alpha+1}}(\mathfrak{M})\models \mathsf{ACA}_0$.

To finish the proof of  $\mathbf{DEF}_{\omega^{\alpha+1}}(\mathfrak{M})\models (\Pi^0_1\mbox{-}\mathsf{CA}_0)_{\omega^{1+\alpha}}$ we need to show that for any $\Pi^0_{\infty}$ formula $A(x,Y)$ whose parameters are from $\mathbf{DEF}_{\omega^{\alpha+1}}(\mathfrak{M})$ there is a set $H\in\mathbf{DEF}_{\omega^{\alpha+1}}(\mathfrak{M})$  such that 
\begin{equation} \label{hier_eq}\mathbf{DEF}_{\omega^{\alpha+1}}(\mathfrak{M})\models(\forall \beta<\omega^{1+\alpha})\;\forall x\;(x\in(H)_{\beta}\mathrel{\leftrightarrow} A(x,(H)_{<\beta})).\end{equation}

Let us fix some $\delta<\omega^{\alpha+1}$ such that all the set parameters of $A$ are sets defined by formulas from $\mathcal{L}_{\delta}$. By (a formalized version of) recursion theorem we define Kalmar elementary sequence  $\langle \varphi_{\beta}(x)\mid \beta\le \omega^{1+\alpha}\rangle$ of G{\"o}del numbers of $\LL_{\omega^{\alpha+1}}$-formulas within $\mathfrak{M}$:
$$\varphi_{\omega\beta}(x)\circeq \exists y (\exists \gamma<\omega\beta)(x=\langle \gamma,y\rangle\land \TT_{\delta+\beta}(\varphi_{\gamma+1}(\underline{x}))),$$ $$\varphi_{\omega\beta+k+1}(x)\circeq\varphi_{\omega\beta+k}(x)\lor \exists y\;(x=\langle \omega\beta+k,y\rangle\land A(y,D_{\varphi_{\omega\beta+k}})).$$

The set $H$ that we want to construct is encoded by the formula $\varphi_{\omega^{1+\alpha}}(x)$. Since the latter is given by a possibly non-standard G{\"o}del number in $\mathfrak{M}$, in order to use it we put it inside an apropriate truth definition. Let $\gamma=\delta+\omega^{\alpha}+1$ and let $\psi(x)$ be the formula $\TT_{\gamma}(\varphi_{\omega^{1+\alpha}}(x))$. Now we just need to show that the set $D_{\psi}$ satisfies the required conditions on $H$, i.e., that
\begin{equation} \label{Def_model_prf_e1}\mathfrak{M}\models \al{ \beta<\omega^{1+\alpha}} \al{x} (x\in (D_{\psi})_{\beta}\mathrel{\leftrightarrow} A(x,(D_{\psi})_{<\beta})).\end{equation}

Since $\mathfrak{M}$ is a model of $\mathsf{R}_{<\omega^{1+\alpha+1}}(\mathsf{IB})$ and $\omega\gamma+\omega<\omega^{1+\alpha+1}$, the model $\mathfrak{M}$ satisfies $\mathsf{R}_{<\omega\gamma+\omega}(\mathsf{R}_{<\omega\gamma+\omega}(\mathsf{IB}))$. Hence, in order to prove (\ref{Def_model_prf_e1}) it is enough to show that within $\mathfrak{M}$, for all (possibly non-standard) $\beta<\omega^{\alpha}$ and all $k$,  there is an $\mathsf{R}_{<\gamma+\omega}(\mathsf{IB})$ proof of the equivalence \begin{equation}\label{Def_model_prf_e2}\forall x\;(x\in (D_{\psi})_{\omega\beta+k}\mathrel{\leftrightarrow} A(x,(D_{\psi})_{<\omega\beta+k})).\end{equation} 

Further, we work within $\mathfrak{M}$. First, notice that $\mathsf{R}_{<\gamma+\omega}(\mathsf{IB})$ proves $$\forall x\in \LL_{\delta+\beta}(\TT_{\delta+\beta}(x)\mathrel{\leftrightarrow}\TT_{\delta+\omega^{\alpha}}(x)).$$ Let us construct an $\mathsf{R}_{<\gamma+\omega}(\mathsf{IB})$ proof of the equivalence $$\forall x\;(x\in (D_{\psi})_{<\omega\beta}\mathrel{\leftrightarrow}\varphi_{\omega\beta}(x)).$$ We achieve this by subsequently proving in $\mathsf{R}_{<\gamma+\omega}(\mathsf{IB})$ the equivalence between the following formulas:
\begin{enumerate}
    \item $x\in(D_{\psi})_{<\omega\beta}$,
    \item $\psi(x)\land \ex{z}\ex{\delta<\omega\beta} x=\langle \delta,z\rangle$,
    \item $\ex{y} \ex{\gamma<\omega^{1+\alpha}}( x=\langle \gamma,y\rangle\land \TT_{\delta+\omega^{\alpha}}(\varphi_{\gamma+1}(\underline{x}))) \land     \ex{z}\ex{\delta<\omega\beta} x=\langle \delta,z\rangle$,
    \item $\ex{y} \ex{\gamma<\omega\beta} (x=\langle \gamma,y\rangle\land \TT_{\delta+\omega^{\alpha}}(\varphi_{\gamma+1}(\underline{x})))$,
    \item $\ex{y} \ex{\gamma<\omega\beta} (x=\langle \gamma,y\rangle\land \TT_{\delta+\beta}(\varphi_{\gamma+1}(\underline{x})))$,
    \item $\varphi_{\omega\beta}(x)$.
\end{enumerate}

Next, by induction on $l$ we produce $\mathsf{R}_{<\gamma+\omega}(\mathsf{IB})$ proofs of $$\forall x\;(x\in (D_{\psi})_{<\omega\beta+l}\mathrel{\leftrightarrow}\varphi_{\omega\beta+l}(x))$$  and of $$\forall x\; (x\in (D_{\psi})_{\omega\beta+l}\mathrel{\leftrightarrow} A(x,D_{\varphi_{\omega\beta+l}})).$$ Therefore,  $\mathsf{R}_{<\gamma+\omega}(\mathsf{IB})$ proves (\ref{Def_model_prf_e2}):
$$\begin{aligned}  x\in (D_{\psi})_{\omega\beta+k}& \mathrel{\leftrightarrow} A(x,D_{\varphi_{\omega\beta+k}})\\
&\mathrel{\leftrightarrow} A(x,(D_{\psi})_{<\omega\beta+k}),
\end{aligned}$$
which concludes the proof of the lemma.\ep

\proof{ of Lemma \ref{ica_ref_ind_w}.}\  We define  an $\omega$-interpretation $\mathcal{D}_{\omega^{\alpha+1}}$ of $(\Pi^0_1\mbox{-}\mathsf{CA})_{\omega^{1+\alpha}}$ in  $\mathsf{R}_{<\omega^{1+\alpha+1}+\omega}(\mathsf{IB})$.  In the interpretation $\mathcal{D}_{\omega^{\alpha+1}}$, sets are interpreted by G{\"o}del numbers of $\LL_{\omega^{\alpha+1}}$-formulas $\varphi(x)$ without other free variables. The interpretation $x\in^{*} \varphi$ of the membership predicate is defined by the formula $\TT_{\omega^{\alpha+1}}(\varphi(\underline{x})).$
In other words, the G{\"o}del number of  $\varphi(x)$ represents the set $D_{\varphi}=\{x\mid \TT_{\omega^{\alpha+1}}(\varphi(\underline{x}))\}$.

$\mathcal{D}_{\omega^{\alpha+1}}$ is a version of the model $\mathbf{DEF}_{\omega^{\alpha+1}}(\mathbb{N})$ formally constructed within $\mathsf{R}_{<\omega^{1+\alpha+1}+\omega}(\mathsf{IB})$. Since in $\mathsf{R}_{<\omega^{1+\alpha+1}+\omega}(\mathsf{IB})$ we can prove the induction for all $\LL_{\omega^{\alpha+1}}$ formulas, as well as a number of natural properties of $\TT_{\omega^{\alpha+1}}$, we can formalize the proof of Lemma \ref{ica_ref_0_w} taking  $\mathcal{D}_{\omega^{\alpha+1}}$ as $\mathfrak{M}$. Thus, we see that $\mathcal{D}_{\omega^{\alpha+1}}$ is an $\omega$-interpretation of $(\Pi^0_1\mbox{-}\mathsf{CA}_0)_{\omega^{1+\alpha}}$. Since the interpretation of any second order formula is a formula of $\mathcal{L}_{\omega^{\alpha+1}+1}$, the full scheme of induction holds under this  interpretation. \ep

\subsection{Interpreting iterated reflection principles in iterated comprehension theories}
In this section we finish the proofs of Theorem \ref{ica_ref_0} and Theorem \ref{ica_ref_ind}. The key technical observation we need is Lemma \ref{pt_from_ca} that will allow us to reason about the translations $\mathcal{E}_X$ within the theories of iterated $\Pi^0_1$-comprehension.

Recall that there is an axiomatization of each $(\Pi_1^0\tCA_0)_{\alpha}$ over $\mathsf{ACA}_0$ by a single sentence asserting that for each $X$ there exists the $\alpha$-th Turing jump $X^{(\alpha)}$. In this assertion we could replace  $\alpha$ with a variable $x$ to obtain a formula $(\Pi_1^0\tCA_0)_{x}$.  This allows us to formulate Lemma \ref{pt_from_ca}.

\begin{lemma} \label{pt_from_ca}$\mathsf{ACA}_0$ proves that
$$\forall X\forall\alpha<\Lambda\;((\Pi^0_1\mbox{-}\mathsf{CA}_0)_{\omega(\alpha+1)}\to \exists! P\; \mathsf{TP}^X_{\alpha}(P)).$$
\end{lemma}

Before proving Lemma \ref{pt_from_ca} we need to `bootstrap' theory $(\Pi^0_1\mbox{-}\mathsf{CA}_0)_{\alpha}$.
Often, in the past works on the theories of iterated comprehension, in addition to the comprehension principles the well-orderness of $\alpha$ was  included as one of the axioms. However, recently D.~Flumini and K.~Sato \cite{FS14} observed that the axiom of iterated $\Pi^0_1$-comprehension along an order $\alpha$ implies the well-foundedness of $\alpha$. We will use a version of their result.
\begin{lemma}\label{SF} In $\mathsf{ACA}_0$ it is provable that  $$\forall \alpha<\Lambda((\Pi^0_1\mbox{-}\mathsf{CA})_{\alpha}\to \mathsf{WO}(\alpha)).$$\end{lemma}
\bp\  Reasoning in $\mathsf{ACA}_0$, we consider some $\alpha<\Lambda$ and under the assumption $(\Pi^0_1\mbox{-}\mathsf{CA})_{\alpha}$ claim that $\alpha$ is a well-ordering. Let us consider any set $X$ of ordinals $< \alpha$ and show that $X$ is either empty or has the least element. 

We claim that there exists $Y\subseteq X$ such that for $\beta\in X$
\begin{equation}\label{viss_yab_par}\beta\in Y\iff \{\gamma \in Y\mid \gamma<\beta\}=\emptyset.\end{equation}
Indeed, we use  $(\Pi^0_1\mbox{-}\mathsf{CA})_{\alpha}$ to construct $Z$ such that, for $\beta<\alpha$,
$$(Z)_{\beta}=\{ x \mid x=0\land (\forall \gamma<\beta) (\gamma\in X\to \langle \gamma,0\rangle\not\in (Z)_{<\beta})\}$$
and  put $$Y=\{\beta<\alpha\mid \beta\in X\mbox{ and }\langle \beta,0\rangle\in Z\}.$$

Notice that by (\ref{viss_yab_par}) any element of $Y$ is its least element, hence $Y$ consists of at most one element. If $Y$ is empty then by (\ref{viss_yab_par})  $$\emptyset=Y=\{\beta\in X \mid \{\gamma \in Y\mid \gamma<\beta\}=\emptyset\}=X$$ and we are done. If $Y=\{\beta\}$ then $\beta$ is the least element of $X$, since (\ref{viss_yab_par}) guarantees that $Y$ contains any element of $X$ that is smaller than $\beta$.\ep

\brem The construction of $Y$ in the proof of Lemma \ref{SF} essentially is a variant of Yablo-Visser paradox \cite{Yab85,Vis89a}. An interesting feature of the proof is that it did not require the use of any kind of fixed points. This contrasts with Visser's construction \cite{Vis89a}, where he used the Diagonal Lemma to show that the paradox is applicable to descending sequences of truth predicates.\erem  

\proof{ of Lemma \ref{pt_from_ca}.}\; We reason in $\mathsf{ACA}_0$. Our goal is to prove 
$$\forall X\forall\alpha<\Lambda\;((\Pi^0_1\mbox{-}\mathsf{CA}_0)_{\omega(\alpha+1)}\to \exists! P\; \mathsf{TP}^X_{\alpha}(P)).$$ We  consider some $\alpha<\Lambda$ and a set $X$. We assume $(\Pi^0_1\mbox{-}\mathsf{CA})_{\omega(\alpha+1)}$ and claim that there exists a unique $P$ such that $\mathsf{TP}^X_{\alpha}(P)$. 

Within this proof it will be useful to consider classes of sentences $C_{\beta,n}\subseteq \LL_{\beta}$  consisting of all $\LL_{\beta}$-sentences of logic depth $n$ . Observe that $\LL_{\alpha}=\bigcup\limits_{\beta\le \alpha,n\in \omega} C_{\beta,n}$.

Let us first prove that, for any two $P,P'$, if $\mathsf{TP}^X_{\alpha}(P)$ and $\mathsf{TP}^X_{\alpha}(P')$ then $P=P'$. For this it is enough to show that $\forall \beta< \alpha \forall n\;( P\cap C_{\beta,n}=P'\cap C_{\beta,n})$. We prove the latter by a straightforward transfinite induction on the values of $\omega\beta+n$. To justify the steps of induction we use the compositionality of both $P$ and $P'$. We could use this kind of argument, since by Lemma \ref{SF} we have $\mathsf{WO}(\omega(\alpha+1))$. 

Now let us construct some $P$ that is a compositional truth definition for $\mathcal{L}_{\alpha}$ over $X$. Using iterated comprehension we define a set $E$ such that for all $\beta\le \alpha$ and $n\in \omega$:
\begin{enumerate}
    \item $(E)_{\omega\beta+n}\subseteq C_{\beta,n}$;
    \item $(E)_{\omega\beta}$ is the set of atomic $\LL_{\beta}$-sentences that are true under the standard interpretation of arithmetical signature, interpretation of the predicate $X(x)$ as the set $X$, and the interpretation of predicates $\TT_{\gamma}$, for $\gamma<\beta$, as the sets $\bigcup \limits_{n\in\omega}(E)_{\omega\gamma+n}$;
    \item each $(E)_{\omega\beta+(n+1)}$ consists of
    \begin{enumerate}
        \item all elements of $(E)_{\omega\beta+n}$;
        \item all sentences $\varphi\land\psi$ such that $\varphi,\psi\in(E)_{\omega\beta+n}$,
        \item all sentences $\lnot\varphi$ such that $\varphi\in C_{\beta,n}$ but $\varphi\not\in(E)_{\omega\beta+n}$,
        \item all sentences $\forall x\;\varphi(x)$ such that $\varphi(\underline{k})\in(E)_{\omega\beta+n}$, for all $k$.
    \end{enumerate}
\end{enumerate} 
We put $P=\bigcup\limits_{\beta\le\alpha,n\in \omega}(E)_{\omega\beta+n}$. With the use of arithmetical transfinite induction over $\omega(\alpha+1)$ it is easy to show that $P$ has the desired properties.\ep

Now we are ready to prove Theorem \ref{ica_ref_0}. For this we will first prove Lemma \ref{IT_into_IC_0} that provides us an embedding of the relevant truth theories into the relevant iterated comprehension theories. Next, we prove Lemma \ref{ica_ref_0_w+} showing that the construction we used in Lemma \ref{ica_ref_0_w} for proving $\LL_{0}$ conservation actually yields the conservation for higher $\LL_{\beta}$'s. Combining Lemma \ref{IT_into_IC_0} and Lemma \ref{ica_ref_0_w+} we immediately obtain Theorem  \ref{ica_ref_0}.

\begin{lemma}\label{IT_into_IC_0} For all $\alpha\geq 0$, $$(\Pi^0_1\mbox{-}\mathsf{CA}_0)_{\omega^{1+\alpha}}\supseteq_{\omega^{1+\alpha+1}} \mathsf{R}_{<\omega^{1+\alpha+1}}(\mathsf{IB}).$$\end{lemma}

\bp\  By Lemma \ref{utb_cons_trn} it is sufficient to check that for any fixed $\beta<\omega^{\alpha+1}$ the theory  $(\Pi^0_1\mbox{-}\mathsf{CA}_0)_{\omega^{\alpha}}$ proves that the $\mathcal{E}_X$-translations of the scheme $\mathsf{UTB}_{ \beta}$ and of the axiom $\TT_{\beta}\mbox{-}\mathsf{RFN}(\mathsf{EA}+\mathsf{UTB}_{<\beta})$ hold for all $X$.

Recall that for any fixed $\beta<\omega^{1+\alpha+1}$ the theory $(\Pi^0_1\mbox{-}\mathsf{CA}_0)_{\omega^{1+\alpha}}$ implies $(\Pi^0_1\mbox{-}\mathsf{CA}_0)_{\beta}$. Thus, by Lemma \ref{pt_from_ca}, for any fixed $\beta<\omega^{\alpha+1}$, $(\Pi^0_1\mbox{-}\mathsf{CA}_0)_{\omega^{\alpha}}$ proves that for any $X$ and any $\beta'\le \beta$ there exist a unique compositional partial truth definition for $\mathcal{L}_{\beta'}$-sentences over $X$. We fix a $\beta<\omega^{\alpha+1}$ and an instance
\begin{equation}\label{UTB_instance}\forall \vec{x}\;(\varphi(\vec{x})\mathrel{\leftrightarrow} \TT_{\beta}(\varphi(\underline{\vec{x}})))\end{equation} of the scheme $\mathsf{UTB}_{ \beta}$. 

Further we reason in $(\Pi^0_1\mbox{-}\mathsf{CA}_0)_{\omega^{1+\alpha}}$. For an arbitrary $X$ we claim that the $\mathcal{E}_X$-translations of $(\ref{UTB_instance})$ and the axiom $\TT_{\beta}\mbox{-}\mathsf{RFN}(\mathsf{EA}+\mathsf{UTB}_{< \beta})$ hold.

There is a unique set $P_{\beta}$ which is a compositional partial truth definition for $\mathcal{L}_{\beta}$-sentences over $X$. Since $P_{\beta}$ is unique, we have $$\forall x\;(x\in P_{\beta}\mathrel{\leftrightarrow}\mathcal{E}_X(\TT_{\beta}(x))).$$ Moreover, for each $\beta'$ occurring in $\varphi$ (such a $\TT_{\beta'}$ is necessarily less than $\beta$), we have $$\forall x\;(x\in P_{\beta}\cap\LL_{\beta'}\mathrel{\leftrightarrow}\mathcal{E}_X(\TT_{\beta'}(x))).$$ Hence, by the compositionality of $P_\gb$ we easily verify (\ref{UTB_instance}).

Due to the compositionality of $P_\gb$ it is easy to check that any logical or non-logical axiom of the $\LL_{\beta}$-theory $\mathsf{EA}+\mathsf{UTB}_{< \beta}$ is an element of $P_\gb$. Thus, given a proof $p$ of an $\LL_{\beta}$ sentence in $\mathsf{EA}+\mathsf{UTB}_{< \beta}$, by a straightforward induction on formulas in $p$ we verify that the conclusion of $p$ is an element of $P_\gb$. Since $P_\gb$ coincides with the  $\mathcal{E}_X$-interpretation of $\TT_{\beta}$, we infer that the $\mathcal{E}_X$-translation of $\TT_{\beta}\mbox{-}\mathsf{RFN}(\mathsf{EA}+\mathsf{UTB}_{< \beta})$ holds.\ep

\begin{lemma}\label{ica_ref_0_w+} For all $\alpha\geq 0$, $$(\Pi^0_1\mbox{-}\mathsf{CA}_0)_{\omega^{1+\alpha}}\subseteq_{\omega^{1+\alpha+1}} \mathsf{R}_{<\omega^{1+\alpha+1}}(\mathsf{IB}).$$\end{lemma}
\bp\  In Lemma \ref{ica_ref_0_w} we proved that for any model $\mathfrak{M}$ of $\mathsf{R}_{<\omega^{1+\alpha+1}}(\mathsf{IB})$, the model $\mathbf{DEF}_{\omega^{\alpha+1}}(\mathfrak{M})$ satisfies $(\Pi^0_1\mbox{-}\mathsf{CA}_0)_{\omega^{1+\alpha}}$. Recall that $D_X$ is the set $\{a\in \mathfrak{M} \mid \mathfrak{M}\models X(a)\}$.  To prove the present lemma we need to verify that for each $\beta<\omega^{\alpha+1}$  the $\mathcal{E}_{D_X}$-interpretation of $\TT_{\beta}$ within $\mathbf{DEF}_{\omega^{\alpha+1}}(\mathfrak{M})$ coincides with $\TT_{\beta}$ within $\mathfrak{M}$. Indeed, this will imply that if $\mathfrak{M}\not\models \varphi$, for some $\LL_{\omega^{\alpha+1}}$-sentence $\varphi$, then $\mathbf{DEF}_{\omega^{\alpha+1}}(\mathfrak{M})\not\models \mathcal{E}_{D_X}(\varphi)$ and hence $\mathbf{DEF}_{\omega^{\alpha+1}}(\mathfrak{M})\not\models \forall Y\; \mathcal{E}_{Y}(\varphi)$. Since the construction works for all $\mathfrak{M}\models \mathsf{R}_{<\omega^{1+\alpha+1}}(\mathsf{IB})$, we will infer the conclusion of the lemma.

Within $\mathbf{DEF}_{\omega^{\alpha+1}}(\mathfrak{M})$ the set $D_{\TT_{\beta}}$ clearly is a compositional truth predicate for $\LL_{\beta}$ over $D_{X}$. Since $\mathbf{DEF}_{\omega^{\alpha+1}}(\mathfrak{M})$ is a model of $(\Pi^0_1\mbox{-}\mathsf{CA}_0)_{\omega(\beta+1)}$, by Lemma \ref{pt_from_ca} the set $D_{\TT_{\beta}}$ is the only compositional truth predicate for $\LL_{\beta}$ over $D_{X}$. Thus, $\TT_{\beta}$ within $\mathfrak{M}$ coincide with the $\mathcal{E}_{D_X}$-interpretation of $\TT_{\beta}$ within $\mathbf{DEF}_{\omega^{\alpha+1}}(\mathfrak{M})$.
\ep

This concludes the proof of Theorem \ref{ica_ref_0}.
Now let us prove Theorem \ref{ica_ref_ind}. We will use the same scheme as in the case of Theorem \ref{ica_ref_0}. First, we prove Lemma \ref{IT_into_IC_ind} that provides an embedding of a theory of iterated truth predicates into a theory of iterated comprehension. Second, in Lemma \ref{ica_ref_ind_w+} we show that the construction of Lemma \ref{ica_ref_ind_w} yields conservation not only for $\LL_0$, but for $\LL_{\omega^{\alpha+1}}$.

\begin{lemma}\label{IT_into_IC_ind}
For all $\alpha\geq 0$, $$(\Pi^0_1\mbox{-}\mathsf{CA})_{\omega^{1+\alpha}}\supseteq_{\omega^{1+\alpha+1}+\omega} \mathsf{R}_{<\omega^{1+\alpha+1}+\omega}(\mathsf{IB}).$$
\end{lemma}
\bp\  It is sufficient to verify in $(\Pi^0_1\mbox{-}\mathsf{CA})_{\omega^{1+\alpha}}$ that the following schemes hold under the $\mathcal{E}_{X}$-translation, for any $X$:
\begin{enumerate}
    \item $\mathsf{UTB}_{\beta}$, for $\beta\le\omega^{\alpha+1}$;
    \item $\LL_{\omega^{\alpha+1}+1}\mbox{-}\mathsf{RFN}(\mathsf{IB})$.
\end{enumerate}

By a straightforward induction we prove  $\forall n\; ((\Pi^0_1\mbox{-}\mathsf{CA}_0)_{\omega^{1+\alpha}n})$ in the theory $(\Pi^0_1\mbox{-}\mathsf{CA})_{\omega^{1+\alpha}}$. Hence, by Lemma \ref{pt_from_ca}, the theory $(\Pi^0_1\mbox{-}\mathsf{CA})_{\omega^{1+\alpha}}$ proves that, for all $\beta<\omega^{\alpha+1}$ and all sets $X$, there exists a unique set $P_{\beta}^X$ such that $\mathsf{TP}^X_{\beta}(P^X_{\beta})$.

Using provable uniqueness of $P^X_{\beta}$, for any particular instance of  $\mathsf{UTB}_{\beta}$ (for some $\beta\le\omega^{\alpha+1}$) we easily prove its   $\mathcal{E}_X$-translation in $(\Pi^0_1\mbox{-}\mathsf{CA})_{\omega^{1+\alpha}}$. Note that there is a difference between the cases of $\beta<\omega^{\alpha+1}$ and $\beta=\omega^{\alpha+1}$. In the latter case we do not have access to a set that is a compositional truth definition for $\mathcal{L}_{\omega^{\alpha+1}}$. The $\mathcal{E}_X$-interpretation of $\TT_{\omega^{\alpha+1}}$ essentially is the union $\bigcup_{\beta<\omega^{\alpha+1}} P^X_{\beta}$ (which is not necessarily a set).

Let us consider an instance of $\LL_{\omega^{\alpha+1}+1}\mbox{-}\mathsf{RFN}(\mathsf{IB})$
\begin{equation}\label{IT_into_IC_ind_e1}\forall \vec{x}\;(\Box_{\mathsf{IB}}\varphi(\underline{\vec{x}})\to \varphi(\vec{x})),\end{equation}
where $\varphi$ is a $\Pi_{n}^{\LL_{\omega^{\alpha+1}+1}}$-formula. Reasoning in $(\Pi^0_1\mbox{-}\mathsf{CA})_{\omega^{1+\alpha}}$, we consider a set $X$ and claim that the $\mathcal{E}_{X}$-translation of (\ref{IT_into_IC_ind_e1}) holds. We assume that for some $\vec{x}$ we have $\Box_{\mathsf{IB}}\varphi(\underline{\vec{x}})$ and prove $\varphi(\vec{x})$. 

By a formalized version of Lemma \ref{utb_cons_trn} we obtain $\Box_{\mathsf{EA}+\mathsf{UTB}_{\le\omega^{\alpha+1}}}\varphi(\underline{\vec{x}})$. Recall that the theory $\mathsf{EA}+\mathsf{UTB}_{\le\omega^{\alpha+1}}$ is $\Pi_2^{\LL_{\omega^{\alpha+1}+1}}$-axiomatizable. Hence, by a standard application of the cut-elimination theorem we get a $\mathsf{EA}+\mathsf{UTB}_{\le\omega^{\alpha+1}}$-proof $p$ of $\varphi(\underline{\vec{x}})$, where all the formulas are from $\Pi_{n+2}^{\LL_{\omega^{\alpha+1}+1}}$.

In order to finish the proof we construct a compositional partial truth definition $\mathsf{Tr}_{n+2}(y)$ for $\Pi_{n+2}^{\LL_{\omega^{\alpha+1}+1}}$-sentences such that $\mathsf{Tr}_{n+2}(\varphi(\underline{\vec{x}}))\mathrel{\leftrightarrow}\mathcal{E}_{X}(\varphi(\vec{x}))$ and the axioms from $\Pi_2^{\LL_{\omega^{\alpha+1}+1}}$-axiomatization of $\mathsf{EA}+\mathsf{UTB}_{\le\omega^{\alpha+1}}$ are true in $\mathsf{Tr}_{n+2}$. Given such a definition, by a straightforward induction on the structure of $p$ we  show that the universal closures of all of them are true, hence $\varphi(\vec{x})$ holds.

To define $\mathsf{Tr}_{n+2}$, we evaluate atomic arithmetical sentences by their standard truth value, we evaluate $X(x)$ by the set $X$, for $\beta<\omega^{\alpha+1}$ we evaluate predicates $\TT_{\beta}$ as $P^X_{\beta}$, and we evaluate $\TT_{\omega^{\alpha+1}}$ as  $\bigcup\limits_{\beta<\omega^{\alpha+1}} P^X_{\beta}$. Next, in a standard manner, we expand $\mathsf{Tr}_{n+2}$ to all $\Pi_{n+2}^{\LL_{\omega^{\alpha+1}+1}}$-formulas (in Appendix \ref{Tr_def_appendix} we construct a partial truth definition for finite languages $\LL$, essentially the same construction works in the case that we are interested in here). Finally, using the compositionality of $P^X_{\beta}$, we verify that all the required axioms are true in $\mathsf{Tr}_{n+2}$.\ep

\begin{lemma}\label{ica_ref_ind_w+} For all $\alpha\geq 0$, $$(\Pi^0_1\mbox{-}\mathsf{CA})_{\omega^{1+\alpha}}\subseteq_{\omega^{1+\alpha+1}+\omega} \mathsf{R}_{<\omega^{1+\alpha+1}+\omega}(\mathsf{IB}).$$\end{lemma}
\bp\  Let us work in $\mathsf{R}_{<\omega^{1+\alpha+1}+\omega}(\mathsf{IB})$. In Lemma \ref{ica_ref_ind_w} we constructed an $\omega$-interpretation $\mathcal{D}_{\omega^{\alpha+1}}$ of $(\Pi^0_1\mbox{-}\mathsf{CA})_{\omega^{1+\alpha+1}}$ in $\mathsf{R}_{<\omega^{1+\alpha+1}+\omega}(\mathsf{IB})$. Let $D_X$ be the set $\{x\mid T_{\omega^{\alpha+1}}(X(x))\}$ within $\mathcal{D}_{\omega^{\alpha+1}}$. To prove the claim we just verify, for any externally given $\beta\le\omega^{\alpha+1}$, that
\begin{equation}\label{ica_ref_ind_w+_e1}\forall x\;(\TT_{\beta}(x)\mathrel{\leftrightarrow} \mathcal{D}_{\omega^{\alpha+1}}(\mathcal{E}_{D_X}(T_{\beta}(x)))).\end{equation} Or in other words, we show that the predicate $\TT_{\beta}$ coincide with  $\mathcal{D}_{\omega^{\alpha+1}}$-interpretation of the $\mathcal{E}_{D_X}$-interpretation of $\TT_{\beta}$.

 Observe that for $\beta<\omega^{\alpha+1}$ the set $D_{\TT_{\beta}}=\{x\mid T_{\omega^{\alpha+1}}(\TT_{\beta}(x))\}$ in $\mathcal{D}_{\omega^{\alpha+1}}$ is  a compositional truth definition for $\LL_{\beta}$ over $D_X$. Hence using the fact that within $\mathcal{D}_{\omega^{\alpha+1}}$, for any $\gamma<\omega^{\alpha+1}$, there exists a unique compositional truth definition for $\LL_{\beta}$ over $D_X$, we easily obtain the equivalence (\ref{ica_ref_ind_w+_e1}) in the case $\beta<\omega^{\alpha+1}$.

To show (\ref{ica_ref_ind_w+_e1}) for the case $\beta=\omega^{\alpha+1}$, 
we prove that, for all $\gamma<\omega^{\alpha+1}$,
\begin{equation}\label{ica_ref_ind_w+_e2}\forall x\in \LL_{\gamma}\;(\TT_{\omega^{\alpha+1}}(x)\mathrel{\leftrightarrow} \mathcal{D}_{\omega^{\alpha+1}}(\mathcal{E}_{D_X}(T_{\omega^{\alpha+1}}(x)))).\end{equation}
Within the $\mathcal{D}_{\omega^{\alpha+1}}$-interpretation, the set $D_{\TT_{\gamma}}$ is the only compositional partial truth definition for $\LL_{\gamma}$ over $X$. Hence,
$$\mathcal{D}_{\omega^{\alpha+1}}(\forall x\in \LL_{\gamma}\;(x\in D_{\TT_{\gamma}}\mathrel{\leftrightarrow} \mathcal{E}_{D_X}(T_{\omega^{\alpha+1}}(x)))),$$
and therefore
\begin{equation}\label{ica_ref_ind_w+_e3}\forall x\in \LL_{\gamma}\;(\TT_{\omega^{\alpha+1}}(\TT_{\gamma}(x))\mathrel{\leftrightarrow} \mathcal{D}_{\omega^{\alpha+1}}(\mathcal{E}_{D_X}(T_{\omega^{\alpha+1}}(x)))).\end{equation}
Using $\LL_{\gamma+1}$-reflection we easily prove that for any $\varphi\in \mathcal{L}_{\gamma}$, we have $\TT_{\omega^{\alpha+1}}(\varphi)\mathrel{\leftrightarrow}\TT_{\omega^{\alpha+1}}(\TT_{\gamma}(\varphi))$. Thus, we derive (\ref{ica_ref_ind_w+_e2}) from (\ref{ica_ref_ind_w+_e3}).\ep
This concludes the proof of Theorem \ref{ica_ref_ind}.

\appendix

\section{Truth definitions for $\Pi^\LL_m$-formulas}\label{Tr_def_appendix}

We consider a signature $\LL$ extending that of arithmetic by finitely many predicate letters. We assume that our language is that of Tait calculus, that is, formulas are obtained from atomic ones and their negations by $\land$, $\lor$, $\al{}$, $\ex{}$.  
Our goal is to show the existence of suitable truth definitions for $\Delta_0^\LL$-formulas. 



\bt \label{Delta_0_Tr} There is a $\Pi_1^\LL$-formula $\Tr$ such that for all $\Delta_0^\LL$-formulas $\varphi(\vec{x})$ 
\benr
\item \label{delta_0_tr_1} $\EA\vdash\forall \vec{x} \;(\Tr(\varphi(\vec{x}))\to\varphi(\vec{x}))$;
\item \label{delta_0_tr_2} $\EA^\LL\vdash\forall \vec{x} \;(\Tr(\varphi(\vec{x}))\mathrel{\leftrightarrow}\varphi(\vec{x}))$.
\eenr 
\et 

Let us define predicate $\Tr$ that is required in Theorem \ref{Delta_0_Tr}. In order to construct $\Tr$ we first introduce a notion of a partial evaluation function $s$ for $\Delta_0^\LL$-sentences.  
Informally, $s$ is called a partial evaluation function if $s$ is a locally correct assignment of truth values to some $\Delta_0$-sentences and numerical values to some closed arithmetical terms. Within $\EA$ we define predicate $\mathsf{EF}(s)$ as the conjunction of the following:
\ben 
\item\label{es_def_1} $s$ is a finite function, whose domain consists of $\Delta_0^\LL$ sentences and closed arithmetical terms;
\item for $\Delta_0$-sentences $\psi\in\dom(s)$, the values $s(\psi)\in\{0,1\}$;
\item for terms $t\in\dom(s)$, the values $s(t)$ are some natural numbers;
\item $0\in\dom(s)\to s(0)=0$;
\item \label{es_def_5}for any closed term $t$, if $(S(t))\in\dom(s)$, then $t\in\dom(s)$ and $s(S(t))=S(s(t))$;
\item for any closed terms $t,v$, if $(t+v)\in\dom(s)$, then $t,v\in\dom(s)$ and $s(t+v)=s(t)+s(v)$;
\item for any closed terms $t,v$, if $(t\times v)\in\dom(s)$, then $t,v\in\dom(s)$ and $s(t\times v)=s(t)s(v)$;
\item \label{es_def_8}for any closed term $t$, if $(\exp(t))\in\dom(s)$, then $t\in\dom(s)$ and $s(\exp(t))=\exp(s(t))$;
\item \label{es_def_9} for any $\psi(x_1,\ldots,x_n)$ that is an $\LL$-predicate or its negation and any closed terms $t_1,\ldots,t_n$, if $\psi(t_1,\ldots,t_n)\in\dom(s)$, then $t_1,\ldots,t_n\in\dom(s)$ and $s(\psi(t_1,\ldots,t_n))=1\eqv \psi(s(t_1),\ldots,s(t_n))$; 
\item for any $\Delta_0^\LL$-sentences $\phi,\psi$, if $(\phi\land\psi)\in\dom(s)$, then $\phi,\psi\in\dom(s)$ and $s(\phi\land\psi)=1\eqv s(\phi)=1\land s(\psi)=1$;
\item for any $\Delta_0^\LL$-sentences $\phi,\psi$, if $(\phi\lor\psi)\in\dom(s)$, then $\phi,\psi\in\dom(s)$ and $s(\phi\lor\psi)=1\eqv s(\phi)=1\lor s(\psi)=1$;
\item for any $\Delta_0^\LL$-formula $\phi(x)$ without other free variables, if $(\forall x\leq t\;\phi)\in\dom(s)$, then $t\in\dom(s)$, $\forall m\le s(t)\;\phi(\num m)\in\dom(s)$, and $s(\al{x\leq t}\phi(x))=1 \eqv \al{m\leq s(t)}s(\phi(\num m))=1$;
\item \label{es_def_13} for any $\Delta_0^\LL$-formula $\phi(x)$ without other free variables, if $(\exists x\leq t\;\phi)\in\dom(s)$, then $t\in\dom(s)$, $\forall m\le s(t)\;\phi(\num m)\in\dom(s)$, and $s(\ex{x\leq t}\phi(x))=1 \eqv \ex{m\leq s(t)}s(\phi(\num m))=1$.
\een
We note that, since $\LL$ is finite, the quatifier over predicate symbols of $\LL$ in clause \ref{es_def_9}. could be formalized in $\LL$ as a finite conjunction. We notice that the external unbounded quantifiers over terms or formulas in clauses \ref{es_def_1}.--\ref{es_def_13}. can, in fact, be bounded by $s$, since the objects they represent are encoded in the finite function $s$. Henceforth $\mathsf{EF}(x)$ could be $\EA$-provably equivalently transformed into a $\Delta_0^\LL$-formula.

Observe that provably in $\mathsf{EA}$, different evaluation function agree with each other on the intersections of their domains.
\begin{lemma}($\mathsf{EA}$)\label{delta_0_tr_agree} If $\mathsf{EF}(s_1)$ and $\mathsf{EF}(s_2)$, then $$s_1\upharpoonright (\dom(s_1)\cap\dom(s_2))=s_2\upharpoonright (\dom(s_1)\cap\dom(s_2)).$$\end{lemma}
\bp\  We reason in $\EA$. Consider some partial evaluation functions $s_1,s_2$. By induction on the length of terms $t\in \dom(s_1)\cap\dom(s_2)$ we prove that $s_1(t)=s_2(t)$. Next, by induction on the number of connective in sentences $\varphi\in \dom(s_1)\cap\dom(s_2)$ we prove that $s_1(\varphi)=s_2(\varphi)$. Notice that despite the proofs of steps of induction rely on the use of additional predicates from $\LL$, the induction formulas are $\Delta_0$. Thus the proofs are formalizable in $\EA$.\ep

Now, in terms of $\mathsf{EF}(x)$, we define
$$\Tr(\phi):= \al{s}(\mathsf{EF}(s) \land \phi\in\dom(s)\to s(\phi)=1).$$

\begin{lemma}\label{delta_0_tr_lm_-1} For each arithmetical term $t(x_1,\ldots,x_n)$
$$\begin{aligned}\EA\vdash \forall s \forall m_1,\ldots,m_n\;(\mathsf{EF}(s)\land &(t(\num m_1,\ldots,\num m_n))\in\dom(s)\to \\ & t(m_1,\ldots,m_n)=s(t(\num m_1,\ldots,\num m_n))).\end{aligned}$$
\end{lemma}
\bp\  By induction on the construction of $t$, using clauses \ref{es_def_5}.--\ref{es_def_8}. of the definition of $\mathsf{ES}$.\ep
\begin{lemma} \label{delta_0_tr_lm}For each $\Delta_0^\LL$ formula $\varphi(x_1,\ldots,x_n)$
$$\begin{aligned}\EA\vdash\forall s\forall m_1,\ldots,m_n\;(\mathsf{EF}&(s)\land (\varphi(\num m_1,\ldots\num m_n)) \in \dom(s) \to\\  & (\varphi(m_1,\ldots,m_n)\eqv s(\varphi(\num m_1,\ldots\num m_n))=1)).\end{aligned}$$
\end{lemma}
\bp\  We prove the lemma by induction on construction of $\varphi$, using clauses \ref{es_def_9}.--\ref{es_def_13}. of the definition of $\mathsf{ES}$. When in the induction step we consider the cases of formulas $\varphi(x_1,\ldots,x_n)$  starting with a bounded quantifier we use Lemma \ref{delta_0_tr_lm_-1}.\ep

\begin{lemma} \label{delta_0_tr_lm2} For each $\Delta_0^\LL$ formula $\varphi(x_1,\ldots,x_n)$ we have
$$\EA^\LL\vdash \forall m_1,\ldots,m_n\exists s\;(\mathsf{EF}(s)\land  \varphi(\num m_1,\ldots,\num m_n)\in s).$$
\end{lemma}
\bp\  First we prove by induction on construction of terms $t(x_1,\ldots,x_n)$ that there is a term $b_t(x_1,\ldots,x_n)$ such that
\begin{equation}\label{delta_0_tr_eq_3}\EA^\LL\vdash \forall m_1,\ldots,m_n\exists s\le b_t(m_1,\ldots,m_n)\;(\mathsf{EF}(s)\land t(\num m_1,\ldots,\num m_n)\in \dom(s)).\end{equation} 

And next, by induction on the construction of $\Delta_0^\LL$ formulas $\varphi(x_1,\ldots,x_n)$, we prove that there is a term $b_\varphi(x_1,\ldots,x_n)$ such that
\begin{equation} \label{delta_0_tr_eq_4}\EA^\LL\vdash \forall m_1,\ldots,m_n\exists s\le b_\varphi(m_1,\ldots,m_n)\;(\mathsf{EF}(s)\land \varphi(\num m_1,\ldots,\num m_n)\in\dom(s)).\end{equation}

Essentially, the induction steps in the proofs of both (\ref{delta_0_tr_eq_3}) and (\ref{delta_0_tr_eq_4}) are reduced to merging partial evaluation functions which are known to exist by the induction assumption. The fact that we could perform this merging could be justified using Lemma \ref{delta_0_tr_agree}. We note that in the proof of induction step the only two cases for $\varphi$, when it is important that we consider the theory $\EA^\LL$ rather than $\EA$, are the cases when the top-most connective in $\varphi$ is a bounded quantifier. The reason for this is that those are the only cases when inside an $\EA^\LL$-argument we need to perform a merging of partial evaluation functions $s_1,\ldots,s_k$, where $k$, from external point of view, might not be a standard number.

From (\ref{delta_0_tr_eq_4}) we immediately derive the lemma.\ep

\proof{ of Theorem \ref{Delta_0_Tr}.} Using Lemma \ref{delta_0_tr_lm} and the definition of $\Tr$ we derive part (\ref{delta_0_tr_1}) of Theorem \ref{Delta_0_Tr}. Combining (\ref{delta_0_tr_1}) with Lemma \ref{delta_0_tr_lm2} we prove part (\ref{delta_0_tr_2}).\ep

The classes $\Pi_n^\LL$ and $\gS_n^\LL$ are defined from $\Delta_0^\LL$ as usual: $\Pi_0^\LL=\Sigma_0^\LL=\Delta_0^\LL$, $\Pi_{n+1}^\LL=\{\al{\vec x}\phi(\vec x):\phi\in\gS^\LL_n\}$, and $\gS_{n+1}^\LL=\{\ex{\vec x}\phi(\vec x):\phi\in\Pi^\LL_n\}$.
The corresponding truth predicates can be inductively defined as follows:
\bi
\item $\Tr(x)$ serves as a truth definition for $\Sigma_0^\LL$; 
\item $\Tr_{\Pi^\LL_{n+1}}(\phi):=\al{\psi\in\Sigma_n^\LL}(\phi\circeq\al{\vec x}\psi(\vec x) \to \al{z}\Tr_{\Sigma_n^\LL}(\psi(\num{(z)_0}, \dots, \num{(z)_{k-1}})));$
\item $\Tr_{\Sigma^\LL_{n+1}}(\phi):=\neg \Tr_{\Pi_{n+1}^\LL}(\neg\phi).$
\ei 
In the last item, $\neg\phi$ is obtained by de Morgan rules from $\phi$, therefore $\phi\in\Sigma_m^\LL$  iff $\neg\phi\in\Pi_{m}^\LL$.
Then we readily obtain the following theorem. 

\bt \label{tr-pi-m}
Let $\Gamma$ be either $\Pi_n^\LL$ or $\Sigma_n^\LL$, $n\geq 1$. For each $\phi(\vec x)\in\Gamma$,
\benr
\item $\EA^\LL\vdash \al{\vec x} (\phi(\vec x)\eqv \Tr_\Gamma(\phi(\num{\vec x})))$;
\item $\EA\vdash \al{\vec x} (\phi(\vec x)\to \Tr_\Gamma(\phi(\num{\vec x})))$ if $\Gamma=\Pi_{2k-1}^\LL$ or $\gS_{2k}^\LL$, for some $k\geq 1$;
\item $\EA\vdash \al{\vec x} (\Tr_\Gamma(\phi(\num{\vec x}))\to \phi(\vec x))$ if $\Gamma=\gS_{2k-1}^\LL$ or $\Pi_{2k}^\LL$, for some $k\geq 1$.
\eenr 
\et 

We are actually using Statement (i) and Statememnt (ii) for $\Gamma=\Pi_1^\LL$.

\bibliographystyle{plain}
\bibliography{ref-all2}

\end{document}